\RequirePackage{fix-cm}
\documentclass{svjour3}                     
\smartqed  
\usepackage{amsmath}
\usepackage{amssymb}
\usepackage{graphicx}
\usepackage{epsfig}
\usepackage{mathrsfs}
\usepackage{epstopdf}
\usepackage{xcolor}
\usepackage{diagbox}
\usepackage{multirow}
\usepackage{booktabs}
\usepackage{amsopn}
\usepackage{longtable}
\usepackage{lipsum}
\usepackage{hyperref}
\usepackage{marvosym}
\usepackage{units}         
\usepackage{mathptmx}      
\usepackage{enumitem}
\usepackage{latexsym}
\usepackage{caption}
\usepackage[noend]{algpseudocode}
\usepackage{algorithmicx,algorithm}
\usepackage{multirow}
\usepackage{amsfonts}
\usepackage{float}
\usepackage{makecell}
\usepackage{bookmark}
\allowdisplaybreaks[4]

\begin{document}

\title{Analyses of the contour integral method for time fractional normal-subdiffusion transport equation
}

\titlerunning{CIM for time fractional normal-subdiffusion transport equation}    



\author{Fugui Ma \textsuperscript{1}\and
        Lijing Zhao  \textsuperscript{2}\and
        Weihua Deng \textsuperscript{1}\and
        Yejuan Wang\textsuperscript{1}
}

\authorrunning{F.~G. Ma et al.} 

\institute{
F.~G. Ma  \at
\email{mafg17@lzu.edu.cn}  \and
L.~J. Zhao (\textrm{\Letter}) \at
\email{zhaolj17@nwpu.edu.cn.} \and
W.~H. Deng \at
\email{dengwh@lzu.edu.cn}\and
Y.~J. Wang \at
\email{wangyj@lzu.edu.cn}\and
\at
 $^{1}$ School of Mathematics and Statistics, Gansu Key Laboratory of Applied Mathematics and Complex Systems, Lanzhou University, Lanzhou 730000, P. R. China.
\at
 $^{2}$ a. School of Mathematics and Statistics, Northwestern Polytechnical University, Xian, 710129, P. R. China. b. Research and Development Institute of Northwestern Polytechnical University in Shenzhen, Shenzhen, 518057, P. R. China.\\
}

\date{Received: date / Accepted: date}

\maketitle

\begin{abstract}
In this work, we theoretically and numerically discuss a class of time fractional normal-subdiffusion transport equation, which depicts a crossover from normal diffusion (as $t\rightarrow 0$) to sub-diffusion (as $t\rightarrow \infty$). Firstly, the well-posedness and regularities of the model are studied by using the bivariate Mittag-Leffler function. Theoretical results show that after introducing the first-order derivative operator, the regularity of the solution can be improved in substance. Then, a numerical scheme with high-precision is developed no matter the initial value is smooth or non-smooth. More specifically, we use the contour integral method (CIM) with parameterized hyperbolic contour to approximate the temporal local and non-local operators, and employ the standard Galerkin finite element method for spatial discretization. Rigorous error estimates show that the proposed numerical scheme has spectral accuracy in time and optimal convergence order in space. Besides, we further improve the algorithm and reduce the computational cost by using the barycentric Lagrange interpolation. Finally, the obtained theoretical results as well as the acceleration algorithm are verified by several 1-D and 2-D numerical experiments, which also show that the numerical scheme developed in this paper is effective and robust.
\keywords{Time fractional equations\and Contour integral method\and Regularity analysis\and Error estimates\and Acceleration algorithm}
\subclass{35R11 \and 35B65 \and 65D30 \and 65M15}
\end{abstract}

\section{Introduction}
\label{sec:Intro}

In recent years, the study of fractional diffusion equations (FDEs) has aroused an upsurge (see cf. \cite{Podlubny99,Metzler00,Metzler04,Kilbas06,Schumer03}, etc.). One class of the most popular FDEs is the time fractional diffusion equations (TFDEs), which are closely related to the time-changing processes. Subodinator is a powerful tool to describe such time-changing process. A subordinator $\mathbb{S}=\{\mathbb{S}_t;t\geq0\}$ with $\mathbb{S}(0)=0$ is a stochastic process in continuous time with non-decreasing paths for which $\mathbb{E}[e^{-\lambda \mathbb{S}_t}]=e^{-t\phi(\lambda)}$, $\lambda>0$ where $\phi$ is a Bernstein function, whose inverse (or hitting time) process is defined as $\mathbb{E}_t=\inf\{s>0;\mathbb{S}_s>t\}$. When $\phi$ is taken as $\lambda^{\beta}$ with $0<\beta<1$, the related subordinator is called the $\beta$-stable subordinator and its corresponding inverse process is called the inverse $\beta$-stable subordinator (see cf. \cite{chen18,Toaldo15}, and the references therein).

Denote $\{\bar{\mathbb{S}}_t;t\geq0\}$ be a drift-less subordinator with Laplace exponent  Bernstein function $\phi_0(\lambda)=\int_{0}^{\infty}(1-e^{-\lambda x})\nu(dx)$ and infinite L\'{e}vy measure $\nu$. Let $\mathbb{E}_t$ be the inverse process or hitting time process of the $\beta$-stable subordinator $\mathbb{S}_t=Kt+\bar{\mathbb{S}}_t$, $K\geq0$ ($K=0$ is drift-less). Then its Laplace exponent Bernstein function $\phi$ has the following representation
\begin{displaymath}
\phi(\lambda)=K\lambda+\int_{0}^{\infty}(1-e^{-\lambda x})\nu(dx).
\end{displaymath}

Let $u(t,r)=\mathbb{E}_r[u_0(\mathbf{B}_{\mathbb{E}_t})]$  be a subordinated stochastic process where $u_0$ is a given initial distribution and $\mathbf{B}(t)$ be Brownian motion on $\mathbb{R}^{d}$ ($d = 1, 2, 3$). Then this stochastic process satisfies the following governing equation\cite{chen18}:
\begin{displaymath}
K\frac{d}{d t}u(t,r)+\partial_{t}^{\beta}u(t,r)=\Delta u(t,r),
\end{displaymath}
which can be widely used to describe many diffusion phenomena in physics, porous medium and hydrology, as well as the reactions or mechanisms in chemistry, engineering, finance and social sciences (see cf. \cite{Metzler00,Metzler04,Schumer03,Zhang14} and the references therein).
Indeed, when $K=0$, the above equation describes a sub-diffusion process; while if $K>0$, it depicts a crossover from normal diffusion (as $t\rightarrow0$) to sub-diffusion diffusion (as $t\rightarrow\infty$). See Appendix A for details.

In this paper, we consider a more general time fractional normal-subdiffusion transport model
\begin{equation}\label{eq:problem}
\left\{\begin{split}
&K\frac{d}{d t}u(t,r)+\partial_{t}^{\beta}u(t,r)+Au(t,r)=f(t,r),
\quad && {\rm in} ~(t,r)\in(0, T]\times\Omega, \\
&u(0,r)=u_0,\quad && r\in\Omega,
\end{split}\right.
\end{equation}
where $0<\beta<1$ denotes the fractional order,  $f(t,\cdot)$: $[0, T]\rightarrow L^{2}(\Omega)$ is a given source term, $\Omega$ is a bounded convex polygonal domain in $\mathbb{R}^{d}$ ($d = 1, 2, 3$) with a boundary $\partial\Omega$, $T>0$ is a given time value, $K\geq0$ is a constant, $A:=-\Delta$: $H_{0}^{1}(\Omega)\cap H^{2}(\Omega)\rightarrow L^{2}(\Omega)$ is the negative Laplacian operator with zero Dirichlet boundary condition,  $\partial_{t}^{\beta}u$ is the fractional Caputo derivative which is given in terms of the Riemann-Liouville one (see cf. \cite{Podlubny99,chen18}, etc.) by
\begin{displaymath}
\partial_{t}^{\beta}u(t,r)=\frac{1}{\Gamma(\beta)}\frac{d}{dt}\int_{0}^{t}(t-s)^{\beta-1}(u(s,r)-u(0,r))ds
\end{displaymath}
with $\Gamma(\cdot)$ being the Euler Gamma function $\Gamma(\lambda):=\int_{0}^{\infty}t^{\lambda-1}e^{-t}dt$.

In the past two decades, a great deal of studies have been focused on the drift-less case $K=0$, which include theoretical analyses and numerical computations.
In theory, the existing works primarily analyze the well-posedness and regularity of the solution (see cf. \cite{Kilbas06,Mu17,Kubica18}, etc.). Numerically, researchers consider on how to effectively approximate the time non-local operators. At present, popular numerical discretization mainly include the following categories: $L_1$-type methods (see cf. \cite{Lin07,Jin16b,Gao14,Deng08,Jiang17}, etc); convolution quadrature (CQ) methods (see cf. \cite{Lubich86,Cuesta06,Zeng15,Deng19}, etc), and other improved versions based on these (see cf. \cite{Yan17,Li21,Stynes17,Mustapha20}, etc.).

For $K>0$,  Problem (\ref{eq:problem}) is expected to improve the modeling delicacy in depicting the anomalous diffusion. In particular, when $K=1$, Problem (\ref{eq:problem}) is the so-called time-fractional mobile-immobile transport model, which can describe the mechanical behavior of anomalous diffusion transport in heterogeneous porous media fine (see \cite{Schumer03,Zhang14}, etc).
In \cite{chen18}, the author analyzes the existence and uniqueness of the solution from the perspective of probability theory by using random representation; based on some given regularity assumptions, \cite{Nikana20} uses the CQ method and the radial basis function-generated finite difference method to solve Problem (\ref{eq:problem}); for $K=1$, \cite{Mustapha20} proposes a second-order accurate $L_1$ scheme over non-uniform time steps,  \cite{Zheng22} develops the averaged $L_1$ compact difference scheme; \cite{Bazhlekova21} studies the relaxation functions for equations with multiple time-derivatives and  discusses the subordination result as well. To the best of our knowledge, the present works on numerical approximation, well-posedness and regularity analysis for Problem (\ref{eq:problem}) are still limited.

\begin{table*}[t!]
\setlength{\abovecaptionskip}{0cm}  
\setlength{\belowcaptionskip}{-0.2cm} 
\scriptsize
\caption{Part summary of convergence rates and regularity assumptions in time for existing time semi-discrete schemes of Problem (\ref{eq:problem}) when $K=0$, where $\tau$ is the time discrete step-spacing.}
\label{Tab:00}
\begin{center}
\begin{tabular}{p{4.0cm}p{4.0cm}c}
  \hline
  \hline
  method & rate  & regularity assumption  \\
  \hline
  \cite{Jin16b} & $\mathcal{O}(\tau)$  & $\forall x\in \Omega$, $u$ is $C^{2}$ in $t$.\\
  \cite{Gao14} & $\mathcal{O}(\tau^{3-\beta})$  & $\forall x\in \Omega$, $u$ is $C^{3}$  in $t$.\\
  \cite{Lin07} & $\mathcal{O}(\tau^{2-\beta})$  & $\forall x\in \Omega$, $u$ is $C^{2}$ in $t$.\\
  \cite{Zeng15}  & $\mathcal{O}(\tau^{2-\beta})$  & $\forall x\in \Omega$, $u$ is $C^{2}$ in $t$.\\
  \cite{Jiang17} & $\mathcal{O}(\tau^{2-\beta})$  & $\forall x\in \Omega$, $u$ is $C^{2}$ in $t$.\\
  \cite{Yan17} & $\mathcal{O}(\tau^{2})$  & $\forall x\in \Omega$, $u$ is $C^{2}$ in $t$.\\
  \cite{Li21} & $\mathcal{O}(\tau^{2})$  & $\forall x\in \Omega$, $u$ is $C^{3}$ in $t$.\\
  \hline\hline
\end{tabular}
\end{center}
\end{table*}

From the works list above, we can see that, compared with theoretical analyses, numerical computations for TFDEs are predominant. Although some high-order schemes have been proposed, strong regularity of solutions should be required at the same time (see Table \ref{Tab:00}), which is inconsistent with the regularity of the physical model itself. Besides, due to the historical dependence of the time non-local operators and its weak singularity at the original point, the storage and computational costs of the traditional numerical schemes are usually expensive. This paper targets on: (i) the study of the well-posedness of Problem (\ref{eq:problem}) for $K>0$ and the regularity of its solution from the perspective of PDE. (ii) high-performance numerical method with low regularity and $\mathcal{O}(1)$ storage requirements in temporal direction. (iii) algorithm acceleration to reduce the computational operations.

To illustrate the time semi-discrete method, we consider the following time-fractional initial-boundary value problem (i.e., Problem (\ref{eq:problem}) with $K=0$)
\begin{equation}\label{eq:iniproblem1}
   \partial_{t}^{\beta}u(t,r) + Au(t,r) = f(t,r)\quad~ {\rm with}~ ~ u(0,r)=u_{0}.
\end{equation}
Based on the fact that the spectrum of $A$ is constrained in a sufficiently small sector, i.e.,
\begin{displaymath}
sp(A)\subset\Sigma_{\theta}:=\left\{z\in\mathbb{C}: z\neq0, |\arg(z)|<\theta\right\}~ ~{\rm with}~\theta\in(0,\pi/2),
\end{displaymath}
and the resolvent $(zI+A)^{-1}: L^{2}(\Omega)\rightarrow L^{2}(\Omega)$ satisfies (see \cite[Theorem 3.7.11]{Arendt11})
\begin{equation}\label{eq:resolvent}
\left\|(zI+A)^{-1}\right\|_{L_2(\Omega)}\leq\frac{C}{1+|z|}\quad{\rm for~ all}~z\in\Sigma_{\theta},~\theta\in(\pi/2,\pi).
\end{equation}
the construction of the time semi-discrete method mainly contains two steps (for convenience, we will ignore the spatial variable identifier in the sequence, i.e., $u(t):=u(t,r)$, $f(t):=f(t,r)$) :
\begin{description}
  \item[\emph{Step I}:]
Under weak restrictions on $u(t)$, there exists a positive constant $\sigma_0$ ($\sigma_0>0$ is often referred to as the convergent abscissa), such that the Laplace transform $\widehat{u}(z)$ exists for ${\rm Re}(z)>\sigma_0$ (i.e., there is $C>0$, $|u(t)|\leq C e^{\sigma_0 t}$)(see \cite{Essah}),
we express the solution $u(t)$ in the form of
\begin{equation}\label{eq:Laplace}
  u(t)=\frac{1}{2\pi i}\int_{\sigma-i\infty}^{\sigma+i\infty}e^{zt}\widehat{u}(z)dz,\quad {\rm Re(z)}>\sigma_{0},
  \end{equation}
and
\begin{displaymath}
\widehat{u}(z)=\left(z^{\beta}I+A\right)^{-1} \left(z^{\beta-1}u_{0}+\widehat{f}(z)\right).
\end{displaymath}
  \item[\emph{Step II}:] Propose an efficient numerical method to approximate the improper integral (\ref{eq:Laplace}). Attracted by the simplicity and efficiency of the contour integral method (CIM) (cf. \cite{Fernandz04,McLean04,Talbot,Sheen00}), we further develop it in this paper. Other numerical methods for the improper integrals can see cf. \cite{Piessens75,Piessens76}, etc.

\end{description}

The basic idea of CIM is, by Cauchy's integral theorem, the original integration path of the inverse Laplace transform (is a vertical line from negative infinity to positive infinity) can be deformed into a contour with ${\rm Re}(z)\rightarrow-\infty$ at each end. Thus, the exponential factor $e^{zt}$ can force a rapid decay of the integrand on the contour, which is greatly beneficial to the fast convergence of the numerical computation to the improper integral (\ref{eq:Laplace}), and avoiding the unfeasible approximation on the vertical line as well as high frequency oscillation.

Let such an appropriate contour be parameterized by
\begin{equation}\label{eq:contour}
\Gamma:z=z(\phi),\ -\infty<\phi<\infty.
\end{equation}
Then the solution (\ref{eq:Laplace}) can be rewritten as
\begin{equation}\label{eq:wcontour}
  u(t)=\frac{1}{2\pi i}\int_{-\infty}^{\infty}e^{z(\phi)t}\widehat{u}(z(\phi))z'(\phi)d\phi.
\end{equation}
Approximating integral (\ref{eq:wcontour}) by mid-point rule with uniform step-spacing $\tau$, we can get
\begin{equation}\label{eq:approx}
 u(t) \approx u^{N}(t)=\frac{\tau}{2\pi i}\sum\limits_{k = -\infty}^{\infty}e^{z(\phi_{k})t}~\widehat{u}(z(\phi_{k}))z'(\phi_{k}),
\end{equation}
where $\phi_{k}=(k+1/2)\tau$ (if trapezoidal rule is used, $\phi_{k}=k\tau$).
Suppose that the contour $\Gamma$ is symmetric with respect to the real axis and Laplace transform of $f(t)$ holds (by Riemann-Schwartz reflection principle, there is $\widehat{f}(\bar{z})=\overline{\widehat{f}(z)}$), then we have $\widehat{u}(\bar{z})=\overline{\widehat{u}(z)}$. After truncating, there is
\begin{equation}\label{eq:finapprox}
  u(t)\approx u^{N}(t) = \frac{\tau}{\pi} {\rm Im} \left\{\sum\limits_{k = 0}^{N-1}e^{z(\phi_{k})t}~\widehat{u}(z(\phi_{k}))z'(\phi_{k})\right\}.
\end{equation}
The discrete scheme (\ref{eq:finapprox}) is known as the CIM of the problem (\ref{eq:iniproblem1}) in the time direction.

To our knowledge, the basic idea of CIM is early described in \cite{Talbot}, where Talbot provides a detailed explanation of the origin and research progress of CIM. Then further discussions on different contour choices as well as fast and accurate integration methods can be found in large literature, such as \cite{Sheen00,Fernandz04,McLean04,Fernandz06,Weideman10,Weideman06,Trefethen00}, etc.
Applications of CIM on FDEs have also been successfully developed. In \cite{Pang16}, CIM with hyperbolic contour is exploited to solve space-fractional diffusion equations. In \cite{Weideman07}, authors discuss CIM with both hyperbolic and parabolic integral contours to solve the sub-diffusion equation. In \cite{McLean10b}, three kinds of methods relating to Laplace transformation and CIM are given to solve a particular TFDE. Recently, CIM with a shifted hyperbolic contour is employed in \cite{Colbrook22b} to analyze the fractional viscoelastic beam equation. In \cite{Li21,Li23}, authors utilize CIM with hyperbolic integral contour to deal with the nonlinear parabolic equations and semilinear subdiffusion equations with rough data, respectively. CIM is also used to solve the time fractional Feynman-Kac equation with two internal states in \cite{Ma23}. Besides, the inverse Laplace transform via contour method has been used to great effect in evaluating (scalar) Mittag-Leffler functions \cite{Garrappa15,McLean21}, which appear ubiquitously when studying TFDEs. Comparing with other traditional numerical methods, CIM, which is spectral accurate and cheap in computation, can belittle the influence of the singularity at the origin, and requires low regularity of the solution. What's more, to compute the solution value at a current moment by CIM, it does not depend on the information on history, and so can perform parallel computing.

The main contributions of this paper are listed as follows:
\begin{itemize}[label=$\bullet$]
  \item Firstly, we discuss the well-posedness and regularities of the new model (\ref{eq:problem}), both in temporal as well as spatial  directions, by use of the bivariate Mittag-Leffler function, and get some novel results (see Subsec.  \ref{subsec:discuss}).

  \item Secondly, we develop a numerical scheme for the time fractional normal-subdiffusion problem  (\ref{eq:problem}) by using the standard Galerkin FEM with continuous piecewise linear functions for space discretization and CIM with parameterized hyperbolic contour for time discretization.

  \item We show the error analysis and optimal convergence estimation in detail. Spatial second-order convergence rate can be reached no matter the initial data is smooth or non-smooth.
 The optimal parameters determining the shape of the contour in CIM are given.

  \item Finally, we develop an acceleration algorithm based on the barycentric Lagrange interpolation approximation, and take several numerical experiments with smooth as well as non-smooth initial data in 1-D and 2-D to demonstrate the efficiency and robustness of our proposed numerical scheme.
\end{itemize}

The rest of this paper is organized as follows: in Sec. \ref{sec:smooth}, the well-posedness and regularity of the solution are discussed; in Sec. \ref{sec:timediscrete} and \ref{sec:Galerkin}, the time and space semi-discrete methods are proposed respectively; in Sec. \ref{sec:fulldiscret}, the fully discrete scheme of Problem (\ref{eq:problem}) is established and the convergence analysis is performed; in Sec. \ref{sec:acceler}, an acceleration algorithm for CIM and its error estimation are presented; in Sec. \ref{sec:Numerical}, numerical experiments are given to demonstrate the theoretical results; some conclusions are made in Sec. \ref{sec:conclusion}.

\section{Smoothness theory}
\label{sec:smooth}
\subsection{Preliminaries}
We begin by introducing some notations about functional spaces that will be adopted in the subsequent context. As we have indicated already in Sec. \ref{sec:Intro}, the operator $A$ is the negative Laplacian with zero Dirichlet boundary condition. Let $(\lambda_{j},\varphi_{j})_{j\in\mathbb{N}}$ be the eigenvalues of $A$ ordered non-decreasingly and the corresponding eigenfunctions normalized in the $L^{2}(\Omega)$ norm. Then the fractional Sobolev space $\dot{H}^{q}(\Omega)$ $(q\geq0)$ is defined as (see cf.\cite{Thomee06}, etc)
\begin{displaymath}
\dot{H}^{q}(\Omega):=\left\{v\in L^{2}(\Omega): \Sigma_{j=1}^{\infty}\lambda^{q}_{j}(v,\varphi_{j})^{2}<\infty\right\}
\end{displaymath}
with the following scalar product and norm
\begin{displaymath}
(u,v)_{\dot{H}^{q}(\Omega)}:=\Sigma_{j=1}^{\infty}\lambda^{q}_{j}(u,\varphi_{j})(v,\varphi_{j}),~ ~  \|v\|_{\dot{H}^{q}(\Omega)}:=(v,v)_{\dot{H}^{q}(\Omega)}^{1/2}~~\forall~u,~v\in\mathrm{span}\{\varphi_{j}\}.\\
\end{displaymath}
It is clear that $\dot{H}^{0}(\Omega)=L^{2}(\Omega)$, $\dot{H}^{1}(\Omega)=H^{1}_{0}(\Omega)$ and  $\dot{H}^{2}(\Omega)=H^{1}_{0}(\Omega)\cap H^{2}(\Omega)$.

We also define the dual space of $\dot{H}^{q}(\Omega)$  as
\begin{displaymath}
\dot{H}^{-q}(\Omega):=\left\{v\in L^{2}(\Omega): \Sigma_{j=1}^{\infty}\lambda^{-q}_{j}(v,\varphi_{j})^{2}<\infty\right\}
\end{displaymath}
with norm
\begin{displaymath}
\|v\|_{\dot{H}^{-q}(\Omega)}:=(v,v)_{\dot{H}^{-q}(\Omega)}^{1/2}~ ~ ~ ~ {\rm and}~ ~ ~ ~ (u,v)_{\dot{H}^{-q}(\Omega)}:=\Sigma_{j=1}^{\infty}\lambda^{-q}_{j}(u,\varphi_{j})(v,\varphi_{j}).
\end{displaymath}

For $q\geq 0$, there is $\dot{H}^{q}(\Omega)\subset L^{2}(\Omega)\subset\dot{H}^{-q}(\Omega)$  (see cf. \cite{Sakamoto11}, etc.).

Besides, for $\theta\in(\pi/2, \pi)$, we define another sector
\begin{displaymath}
\Sigma_{\theta,\delta}:=\left\{z\in\mathbb{C}: |z|\geq\delta>0, |\arg(z)|< \theta\right\},
\end{displaymath}
and an integral contour $\Gamma_{\theta,\delta}\subset\mathbb{C}$
\begin{displaymath}
\Gamma_{\theta,\delta}:=\{z\in\mathbb{C}:|z|=\delta,|\arg(z)|\leq \theta\}\cup\left\{z\in\mathbb{C}:z=re^{\pm i\theta}, 0<\delta\leq r<\infty\right\},
\end{displaymath}
which oriented with an increasing imaginary part and $i^2=-1$.

Next, we perform the well-posedness and regularity analysis of Problem (\ref{eq:problem}). Before this, we remark that throughout this paper $c$ and $C$ denote positive constants, not necessarily the same at different occurrences, which are independent of the functions involved.

\subsection{Representation of solution}\label{subsection2.2}
\label{subsec:regularity}
Let $u(t)$ be the solution of Problem (\ref{eq:problem}). By taking the Laplace transform on the equation in (\ref{eq:problem}) and perform simple operations \cite{Podlubny99}, we get
\begin{equation}\label{eq:equality}
\left(\left(Kz+z^{\beta}\right)I+A\right)\widehat{u}(z)=\left(K+z^{\beta-1}\right)u_0+\widehat{f}(z).
\end{equation}
Denote
\begin{equation}\label{eq:denotew}
\eta(z):=Kz+z^{\beta}~~{\rm and}~~W(z):=\frac{\eta(z)}{z}(\eta(z)I+A)^{-1}.
\end{equation}
Then  Eq. (\ref{eq:equality}) can be rewritten as
\begin{equation}\label{eq:integrand}
\widehat{u}(z)=W(z)u_0+z\eta^{-1}(z)W(z)\widehat{f}(z).
\end{equation}
Transforming the inverse Laplace transform, the solution $u(t)$ can be expressed as
\begin{equation}\label{eq:integral}
u(t)=\frac{1}{2\pi i}\int_{\sigma_{0}-i\infty}^{{\sigma_{0}+i\infty}}e^{zt}\widehat{u}(z)dz,~ ~ {\rm Re}(z)>\sigma_{0}.
\end{equation}
From Remark \ref{thm:remark1} below, the integrand $\widehat{u}(z)$ in Eq. (\ref{eq:integral}) is analytical for $z\in\Sigma_{\theta,\delta}$ and any $\delta>0$.
Therefore, we can deform the integral contour from $\sigma_{0}+i\mathbb{R}$ to $\Gamma_{\theta,\delta}$, which reads
\begin{equation}\label{eq:integra_u_G}
u(t)=\frac{1}{2\pi i}\int_{\Gamma_{\theta,\delta}}e^{zt}\widehat{u}(z)dz~ ~ {\rm with}~\theta\in(\pi/2,\pi).  \\
\end{equation}
We formally name the above solution as the mild solution of Problem (\ref{eq:problem}).
\begin{remark}\label{thm:remark1}
The integrand function $\widehat{u}(z)$ in Eq. (\ref{eq:integral}) is analytical for $z\in\Sigma_{\theta}$ with $\theta\in(\pi/2,\pi)$. In fact, according to the expression of $\widehat{u}(z)$ in Eq. (\ref{eq:integrand}), there are two possible cases that may cause singularity: one is $z^{\beta-1}$; another is the operator $(\eta(z)I+A)^{-1}$. For the previous one, obviously, $z=0$ is a singularity point. As for $(\eta(z)I+A)^{-1}$, singularity occurs only if  ${\rm Re}(\eta(z))=-\lambda_j<0$ and ${\rm Im}(\eta(z))=0$ are satisfied at the same time. But, it is impossible. Because, by letting $z=re^{i\zeta}$, $r\geq0$ and $\zeta\in[0,\pi]$, it can be seen that when $r=0$ or $\zeta=0$, there is ${\rm Im}(\eta(z))= Kr\sin(\zeta)+r^{\beta}\sin(\beta\zeta)=0$, while ${\rm Re}(\eta(z))=Kr\cos(\zeta)+r^{\beta}\cos(\beta\zeta)\geq0$, which contraries to ${\rm Re}~\eta(z)<0$. Thus, $\widehat{u}(z)$ is only singular at $z=0$ and analytic for $z\in\Sigma_{\theta}$, with $\theta\in(\pi/2,\pi)$.
\end{remark}
\begin{proposition}\label{thm:prop1}
Let $\eta(z)$ and $W(z)$ be defined in Eq. (\ref{eq:denotew}).
Let $\delta>\max\left\{1,(\frac{2}{K})^{\frac{1}{1-\beta}}\right\}$ be large enough. For $K>0$, we have the following assertions,
\begin{description}
  \item[(I)]for fixed $\theta'\in(\frac{2\pi}{3},\pi)$, there are $\eta(z)\in\Sigma_{\theta}$ with $\theta\in(\frac{\pi}{2},\pi)$ and
\begin{equation}\label{ieq:inequality1}
c|z|<|\eta(z)|\leq C |z|~ ~ {\rm for~ all}~ z\in\Sigma_{\theta',\delta}.
\end{equation}
  \item[(II)] the operator $W(z)$: $L^{2}(\Omega)\rightarrow L^{2}(\Omega)$ is well-defined and bounded for $z\in\Sigma_{\theta',\delta}$, satisfying
\begin{equation}\label{ieq:inequality2}
\left\|W(z)\right\|_{L^{2}(\Omega)}\leq C|z|^{-1}~ {\rm and}~ \left\|AW(z)\right\|_{L^{2}(\Omega)}\leq C.
\end{equation}
\end{description}
\end{proposition}

\begin{proof}
For (I), according to the condition that $\delta$ satisfies, for all $z\in\Sigma_{\theta',\delta}$, there holds
\begin{equation}\label{fig111}
    \sin\left(\left|\arg(Kz+z^{\beta})-\arg(Kz)\right|\right)\cdot K|z|
 \leq\left|z\right|^{\beta},
\end{equation}
the geometric interpretation of which is diagrammed by Fig. \ref{Eq:proppp1}.
\begin{figure}[h]
  \centering
  \includegraphics[height=3.5cm,width=8.8cm]{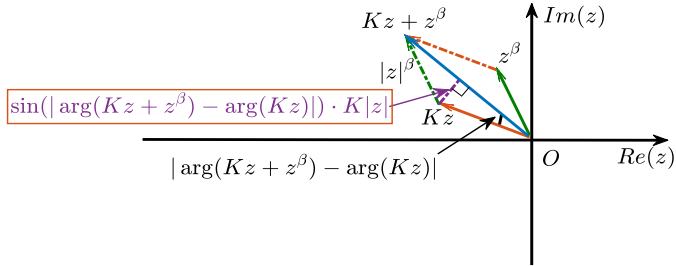}\\
  \caption{Diagram to the inequality (\ref{fig111}).}
  \label{Eq:proppp1}
\end{figure}
Therefore,
\begin{displaymath}
\left|\arg\left(Kz+z^{\beta}\right)-\arg(Kz)\right|
\leq \arcsin\left(\frac{1}{K|z|^{1-\beta}}\right)
\leq \arcsin\left(\frac{1}{K\delta^{1-\beta}}\right)<\frac{\pi}{6}.
\end{displaymath}
Consequently, if $z\in\Sigma_{\theta'}$, $\theta'\in(\frac{2\pi}{3},\pi)$, then $\eta(z)$ belongs to $\Sigma_{\theta}$ with $\theta\in(\frac{\pi}{2},\pi)$. In addition, for all $z\in\Sigma_{\theta',\delta}$, we have
\begin{displaymath}
\frac{|Kz+z^{\beta}|}{|z|}\leq \frac{K|z|+|z|^{\beta}}{|z|}\leq K+\frac{1}{|z|^{1-\beta}}\leq C
\end{displaymath}
and
\begin{displaymath}
\frac{|Kz+z^{\beta}|}{|z|}\geq \frac{\left|K|z|-|z|^{\beta}\right|}{|z|}\geq \frac{K}{2}.
\end{displaymath}
Therefore, we obtain that $c|z|\leq|\eta(z)|\leq C|z|$. The proof of the Proposition \ref{thm:prop1}(I) is completed.

For (II), by Proposition \ref{thm:prop1} (I), for $z\in\Sigma_{\theta',\delta}$, there holds $\eta(z)\in\Sigma_{\theta}$. Therefore, the resolvent estimation of operator $(\eta(z)I+A)^{-1}$ (see Eq. (\ref{eq:resolvent})) satisfies
\begin{displaymath}
\left\|(\eta(z)I+A)^{-1}\right\|_{L^{2}(\Omega)}\leq \frac{C}{1+|\eta(z)|}.
\end{displaymath}
Let
\begin{displaymath}
z\eta^{-1}(z)(\eta(z)I+A)\phi=\psi,
\end{displaymath}
which can be reformulated as
\begin{displaymath}
(\eta(z)I+A)\phi=z^{-1}\eta(z)\psi,
\end{displaymath}
and thus
\begin{displaymath}
\phi=z^{-1}\eta(z)(\eta(z)I+A)^{-1}\psi.
\end{displaymath}
By the resolvent estimation mentioned above, we then have
\begin{displaymath}
\|\phi\|_{L^{2}(\Omega)}\leq C|z|^{-1}\|\psi\|_{L^{2}(\Omega)}.
\end{displaymath}
Now we can deduce that the operator $W(z):$ $L^{2}(\Omega)\rightarrow L^{2}(\Omega)$ is well-defined and bounded.
Since $A(\eta(z)I+A)^{-1}=I-\eta(z)(\eta(z)I+A)^{-1}$, by Proposition \ref{thm:prop1} (I), there also holds
\begin{displaymath}
\|AW(z)\|_{L^{2}(\Omega)}
=\left\|z^{-1}\eta(z)-z^{-1}\eta(z)\eta(z)(\eta(z)I+A)^{-1}\right\|_{L^{2}(\Omega)}\leq C.
\end{displaymath}
The proof of this proposition is completed.  \qed
\end{proof}

Here, we remark that unless otherwise specified, when $z\in\Sigma_{\theta,\delta}$ in the following analyses,  $\theta$ is always selected as $\theta\in(\frac{2\pi}{3},\pi)$ and $\delta$ satisfies the condition mentioned in Proposition \ref{thm:prop1}.

\subsection{Solution theory}
In this section, we shall analyze  the well-posedness and regularity of Problem (\ref{eq:problem}) with $K>0$. For the case of $K=0$, one can refer to the literature \cite{Sakamoto11,Jin18,Mu17}, etc..

\subsubsection{Well-posedness and regularity}
Since $\{(\lambda_{j},\varphi_{j})\}_{j\in\mathbb{N}}$ are the eigenvalues ordered non-decreasingly and the normalized eigenfunctions of $A$, there holds
\begin{equation}
-\Delta\varphi_{j}=\lambda_{j}\varphi_{j}~ {\rm in}~ \Omega,~ {\rm and}~ \varphi_{j}=0~ {\rm on}~ \partial\Omega,~ {\rm for~ all}~ j\in\mathbb{N}_{+}.
\end{equation}
By using the standard separation of variables and eigenvalue expansions, Problem (\ref{eq:problem}) reduces to the following initial value problem
\begin{equation}\label{pro:ODE}
K\frac{\partial}{\partial t}u_j(t)+\partial_{t}^{\beta}u_j(t)+\lambda_ju_j(t)=(f(t),\varphi_{j})~ ~  {\rm with}~ ~ u_{j,0}=(u_0,\varphi_{j}), ~ {\rm for~ all}~ j\in\mathbb{N}_{+}.
\end{equation}
Taking Laplace transform for the above problem, and combining
\begin{displaymath}
\frac{K+z^{\beta-1}}{Kz+z^{\beta}+\lambda_j}=\frac{z^{-1}}{1+\frac{1}{K}z^{\beta-1}+\frac{\lambda_j}{K}z^{-1}}
+\frac{1}{K}\cdot\frac{z^{\beta-2}}{1+\frac{1}{K}z^{\beta-1}+\frac{\lambda_j}{K}z^{-1}},
\end{displaymath}
with the operational method in fractional calculus mentioned in Ref. \cite{Luchko99}, the solution $u(t)$ of problem (\ref{eq:problem}) can be formally represented by
\begin{equation}\label{eq:solutions}
u(t)=E(t)u_{0}+\int_0^{t}\overline{E}(t-\tau)f(\tau)d\tau,
\end{equation}
where
\begin{equation}
\overline{E}(t)\chi:=\sum\limits_{j=1}^{\infty}E^1_{(1-\beta,1),1}\left(\frac{-1}{K}t^{1-\beta},\frac{-\lambda_j}{K}t\right)(\chi,\varphi_j)\varphi_j,
\end{equation}
\begin{equation}\begin{split}
E(t)\chi:=&\sum\limits_{j=1}^{\infty}\left[E^1_{(1-\beta,1),1}\left(\frac{-1}{K}t^{1-\beta},\frac{-\lambda_j}{K}t\right) \right.\\
&\left.+\frac{1}{K}t^{1-\beta}E^1_{(1-\beta,1),2-\beta}\left(\frac{-1}{K}t^{1-\beta},\frac{-\lambda_j}{K}t\right)\right](\chi,\varphi_j)\varphi_j,
\end{split}\end{equation}
and $E^1_{(\alpha,\beta),\gamma}(z_1,z_2)$ is the bivariate Mittag-Leffler function defined in Eq. (\ref{def:BIML}), whose properties are shown in the lemmata below (see Appendix B.2 for the proof of Lemma \ref{Lem:MandL}).
\begin{lemma}\label{Lem:MandL}
Let $\omega_1$, $\omega_2$, $\alpha$, $\beta$, $\gamma\in \mathbb{R}$ with $0<\alpha<\beta\leq 1$, $\gamma\geq0$ and $\omega_1$, $\omega_2<0$. For $t\in[0,\infty)$, there exists a positive constant $C$, such that
\begin{displaymath}
\left|E^1_{(\alpha,\beta),\gamma}\left(\omega_1 t^{\alpha}, \omega_2 t^{\beta}\right)\right|\leq\frac{C}{1+|\omega_2 t^{\beta}|}.
\end{displaymath}
\end{lemma}
\begin{lemma}[cf. \cite{Fernandez20,Bazhlekova21}] \label{Lemma:ML}
Let $\alpha$, $\beta$, $\gamma$, $\omega_1$, $\omega_2$ $\in \mathbb{C}$ with ${\rm Re(\alpha)}$, ${\rm Re(\beta)}$, ${\rm Re(\gamma)}>0$. Then the relationship between fractional calculus of $E_{(\alpha,\beta),\gamma}(\omega_1 t^{\alpha}, \omega_2 t^{\beta})$ and the bivariate Mittag-Leffler function is
\begin{displaymath}\begin{split}
D^{\mu_1}_{t}\left[t^{\gamma-1}E^1_{(\alpha,\beta),\gamma}\left(\omega_1 t^{\alpha}, \omega_2 t^{\beta}\right)\right]=t^{\gamma-\mu_1-1}E^1_{(\alpha,\beta),\gamma-\mu_1}\left(\omega_1 t^{\alpha},\omega_2 t^{\beta}\right),
\end{split}\end{displaymath}
for $\mu_1\in\mathbb{C}$ (fractional integral if ${\rm Re(\mu_1)}<0$, or fractional derivative if ${\rm Re(\mu_1)}\geq0$), where $D^{\mu_1}_{t}$ denotes the Riemann-Liouville fractional derivative (integral) defined as $D^{\mu_1}_{t}u(t)=\frac{1}{\Gamma(\mu_1)}\frac{d^k}{dt^k}\int_{0}^{t}(t-\tau)^{\mu_1-1}u(\tau)d\tau$ with $k=\lfloor {\rm Re(\mu_1)}\rfloor+1$.
\end{lemma}
\begin{remark}
We note that, by the definitions of the two fractional operators $D^{a}_{t}$ and $\partial^{a}_{t}$, it is easy to verify when $\gamma\geq 1$ the equation in Lemma \ref{Lemma:ML} still holds for the case of Caputo fractional derivative.
\end{remark}
\begin{definition}[cf. \cite{Sakamoto11}]\label{def:definition}
We call $u(\cdot,t)$ a weak solution to Problem (\ref{eq:problem}) with zero Dirichlet boundary condition if (\ref{eq:problem}) holds in $L^{2}(\Omega)$ and $u(\cdot,t)\in H_{0}^{1}(\Omega)$ for almost all $t\in(0,T)$ and $u(\cdot,t)\in C([0,T],\dot{H}^{-p}(\Omega))$, $\lim_{t\rightarrow0}\|u(\cdot,t)-u_0\|_{\dot{H}^{-p}(\Omega)}=0$  with $p>0$.
\end{definition}

Now we are ready to give the well-posedness and regularity of the solution to Problem (\ref{eq:problem}).
\begin{theorem}\label{thm:regularityhomo}
Let $0<\beta<1$, $K>0$ and $f\equiv0$.
\begin{description}
  \item[(I)] If $u_0\in L^{2}(\Omega)$, then the weak solution $u(t)$ to Problem (\ref{eq:problem}) given in Eq. (\ref{eq:solutions}) is unique, satisfying $u(t)\in C^{1}((0,T]; L^{2}(\Omega))\cap C([0,T]; {L}^{2}(\Omega))\cap C((0,T]; H^{2}(\Omega)\cap H^{1}_{0}(\Omega))$, with $\partial_t^{\beta}u(t)\in C((0,T]; L^{2}(\Omega))$, and
        \begin{equation}\left\{\begin{split}\label{eq:nonsmoothest}
          &\left\|u(t)\right\|_{C\left([0,T];L^{2}(\Omega)\right)}\leq c \|u_0\|_{L^{2}(\Omega)},\\
          &\|u(t)\|_{\dot{H}^{2}(\Omega)}+\left\|\partial_t u(t)\right\|_{L^{2}(\Omega)}
               + t^{\beta-1}\left\|\partial_t^{\beta}u(t)\right\|_{L^{2}(\Omega)}
          \leq c t^{-1}\|u_0\|_{L^{2}(\Omega)}, t>0.
        \end{split}\right.
        \end{equation}
  \item[(II)] If $u_0\in H^{1}_{0}(\Omega)$, then the unique weak solution $u(t)$  given in  Eq. (\ref{eq:solutions}) belongs to \\ $C^{1}((0,T]; L^{2}(\Omega))\cap L^{2}(0,T; H^{2-\varepsilon}(\Omega)\cap H_{0}^{1}(\Omega))\cap C((0,T]; H^{2}(\Omega)\cap H^{1}_{0}(\Omega))$ for any $\varepsilon>0$, with $\partial^{\beta}_t u(t)\in C((0,T]; L^{2}(\Omega))\cap L^{2}((0,T)\times\Omega)$, and there exists a positive constant $c>0$ such that
      \begin{equation}\left\{\begin{split}\label{Hsmoothest}
      &\|u(t)\|_{L^{2}\left(0,T; H^{2-\varepsilon}(\Omega)\cap H_{0}^{1}(\Omega)\right)}\leq c\|u_0\|_{\dot{H}^{1}(\Omega)},\\
      &\|u(t)\|_{C\left((0,T]; \dot{H}^{2}(\Omega)\right)}
          +\left\|\partial_t u(t)\right\|_{C\left((0,T];L^{2}(\Omega)\right)}
          +t^{\beta-1}\left\|\partial_t^{\beta}u(t)\right\|_{C\left((0,T];L^{2}(\Omega)\right)}\\
      &\leq ct^{-1/2}\|u_0\|_{\dot{H}^{1}(\Omega)}.
      \end{split}\right.
      \end{equation}
  \item[(III)] Furthermore, if $u_0\in H^{2}(\Omega)\cap H^{1}_{0}(\Omega)$, then the unique weak solution $u(t)$  given in  Eq. (\ref{eq:solutions}) belongs to $C^{1}([0,T];L^{2}(\Omega))\cap C([0,T]; H^{2}(\Omega)\cap H^{1}_{0}(\Omega))$, with $\partial^{\beta}_tu\in C([0,T]; L^{2}(\Omega))$, and following  priori estimate holds:
        \begin{equation}\begin{split}\label{est:smooth}
          &\|u(t)\|_{C\left([0,T],\dot{H}^{2}(\Omega)\right)}
          +\left\|\partial_t u(t)\right\|_{C\left([0,T];L^{2}(\Omega)\right)}
          +\left\|\partial_t^{\beta}u(t)\right\|_{C\left([0,T];L^{2}(\Omega)\right)}
          \leq c \|u_0\|_{\dot{H}^{2}(\Omega)}.
        \end{split}\end{equation}
\end{description}
\end{theorem}

\begin{proof}
We firstly note that (\ref{eq:nonsmoothest}) and (\ref{Hsmoothest}), which  extend the resluts of Theorem 2.1 in \cite{Li15}, can also be obtained by a similar way. However, those techniques are no more suitable to get (\ref{est:smooth}) since we can not obtain any estimation about $\frac{d}{dt}u(\cdot,t)$ or $\partial^{\beta}u(\cdot,t)$ by applying them directly. So, here we only show the proof of (III).

In fact, from Theorem \ref{thm:regularityhomo} (I), we have known that the mild solution $u(\cdot,t)$ defined in Eq. (\ref{eq:integral}) is indeed a weak solution to Problem (\ref{eq:problem}).
Note that $\eta(z)A^{-1}(\eta(z)+A)^{-1}=A^{-1}-(\eta(z)+A)^{-1}$ and $\int_{\Gamma_{\theta, \delta}}e^{zt}z^{m-1}dz=0$ with $m\geq1$. Thus, we have
\begin{align*}
\frac{d}{d t}u(t)
=\frac{1}{2\pi i}\int_{\Gamma_{\theta, \delta}}e^{zt}\eta(z)A^{-1}(\eta(z)I+A)^{-1}Au_0dz
=-\frac{1}{2\pi i}\int_{\Gamma_{\theta, \delta}}e^{zt}(\eta(z)I+A)^{-1}Au_0dz.
\end{align*}
Then by Proposition \ref{thm:prop1}, we get
\begin{align*}
\left\|\frac{d}{d t}u(t)\right\|_{L^{2}(\Omega)}
&\leq \frac{c}{2\pi}\int_{\Gamma_{\theta, \delta}}e^{|z|t\cos(\theta)}\left\|(\eta(z)I+A)^{-1}\right\|_{L^{2}(\Omega)}\left\|Au_0\right\|_{L^{2}(\Omega)}|dz|\\
&\leq \frac{c}{2\pi}\int_{\Gamma_{\theta, \delta}}e^{|z|t\cos(\theta)}|z|^{-1}\left\|Au_0\right\|_{L^{2}(\Omega)}|dz|
\leq c \left\|u_0\right\|_{\dot{H}^{2}(\Omega)}.
\end{align*}
Therefore we deduce $u(t)\in C^{1}([0,T];L^{2}(\Omega))\cap C([0,T];H^{2}(\Omega)\cap H^{1}_{0}(\Omega))$. Since $f\equiv 0$, we can see that $\left\|\partial_t^{\beta}u(t)\right\|_{L^{2}(\Omega)}\leq C\|Au(t)\|_{L^2(\Omega)}+C\left\|\frac{d}{d t}u(t)\right\|_{L^{2}(\Omega)}
\leq C \left\|u_0\right\|_{\dot{H}^{2}(\Omega)}$, and $\partial_t^{\beta}u(t)\in C([0,T];L^{2}(\Omega))$. Finally, the estimation in Eq. (\ref{est:smooth}) is obtained. \qed
\end{proof}

Next, we turn to the inhomogeneous case with vanishing initial value.
\begin{theorem}\label{est:reguf}
Let $0<\beta<1$, $K>0$ and $u_0\equiv0$. If $f\in L^{\infty}(0,T;L^{2}(\Omega))$, then $u(t)$ defined in Eq. (\ref{eq:solutions}) is the unique weak solution to Problem (\ref{eq:problem}), belonging to $C^{1}([0,T];L^{2}(\Omega))\cap L^2(0,T;H^{2}\cap H_0^{1}(\Omega))$, $\partial_t^\beta u(t)\in L^2((0,T)\times \Omega)$ with $\partial_tu(t)\in L^{2}((0,T)\times \Omega)$. In particular, for $p>d/2$ ($d=1,2$ or $3$ is the spatial dimension), $u(t)\in C([0,T];\dot{H}^{-p}(\Omega))$, and
\begin{displaymath}
\lim\limits_{t\rightarrow0}\|u(t)-u_0\|_{\dot{H}^{-p}(\Omega)}=0.
\end{displaymath}
Furthermore, there exists a positive constant $C$ such that following estimate holds:
\begin{equation}\label{sourceterm}
\left\|u(t)\right\|_{L^{2}\left(0,T;\dot{H}^{2}(\Omega)\right)}+\left\|\partial_tu(t)\right\|_{L^{2}((0,T)\times\Omega)}+\left\|\partial_t^{\beta}u(t)\right\|_{L^{2}((0,T)\times\Omega)}
\leq c\|f\|_{L^{2}\left((0,T)\times\Omega\right)}.
\end{equation}
\end{theorem}

\begin{proof}
For $u_0\equiv0$, according to the expression of $u(t)$ in Eq. (\ref{eq:solutions}), we have
\begin{equation}\label{inhomosolution}
u(t)=\sum\limits_{j=1}^{\infty}\int_{0}^{t}E^1_{(1-\beta,1),1}\left(\frac{-1}{K}(t-\tau)^{1-\beta},
     \frac{-\lambda_j}{K}(t-\tau)\right)\left(f(\tau),\varphi_{j}\right)d\tau\cdot\varphi_{j}.
\end{equation}
By Lemma \ref{Lem:MandL} and Young's inequality, for any $\nu\in(0,1)$, there holds
\begin{equation}\begin{split}\label{ieq:Auinhomo}
\left\|u(t)\right\|^{2}_{\dot{H}^{2}(\Omega)}
=&\sum\limits_{j=1}^{\infty}\lambda_j^{2}\left(\int_{0}^{t}E^1_{(1-\beta,1),1}\left(\frac{-1}{K}(t-\tau)^{1-\beta},
   \frac{-\lambda_j}{K}(t-\tau)\right)\left(f(\tau),\varphi_{j}\right)d\tau\right)^{2}                           \\
\leq&C\sum\limits_{j=1}^{\infty}\left(\int_{0}^{t}\frac{\lambda_j(t-\tau)^{\nu}(t-\tau)^{-\nu}}
      {\left(1+\frac{\lambda_j(t-\tau)}{K}\right)}\left|\left(f(\tau),\varphi_{j}\right)\right|d\tau\right)^{2}  \\
\leq&C\sum\limits_{j=1}^{\infty}\left(\int_{0}^{t}(t-\tau)^{-\nu}d\tau\right)^2
     \left(\int_{0}^{t}\left|\left(f(\tau),\varphi_{j}\right)\right|^2d\tau\right)                               \\
\leq&Ct^{2-2\nu}\int_{0}^{T}\left\|f(t)\right\|^2_{L^{2}(\Omega)}dt
\leq C\left\|f(t)\right\|^2_{L^{2}((0,T)\times \Omega)},
\end{split}\end{equation}
where the second inequality holds because  $\sup_{j\in\mathbb{N}}\frac{\lambda_j(t-\tau)^{\nu}}{1+\frac{\lambda_j(t-\tau)}{K}}\leq C'$ uniformly for $\nu\in(0, 1)$. Therefore, we deduce that $Au(t)\in C([0,T];L^{2}(\Omega))$, which implies $u(t)\in C([0,T];H^{2}(\Omega)\cap H_{0}^{1}(\Omega))$.
Furthermore, we can see that $\|u(t)\|_{L^2(0,T;\dot{H}^{2}(\Omega))}\leq C\|f\|_{L^{2}\left((0,T)\times\Omega\right)}$, i.e., $u(t)\in L^2(0,T;H^{2}\cap H_0^{1}(\Omega))$.

Next, according to the differential properties of the convolution (i.e., $\partial_t [f_1(t)\ast f_2(t)]=\partial_tf_1(t)\ast f_2(t)+f_1(0)f_2(t)$), there holds
\begin{equation}\begin{split}\label{eq:AAa}
 &\frac{d}{dt}\int_{0}^{t}E^1_{(1-\beta,1),1}\left(\frac{-1}{K}(t-\tau)^{1-\beta},
   \frac{-\lambda_j}{K}(t-\tau)\right)\left(f(\tau),\varphi_{j}\right)d\tau             \\
=&\int_{0}^{t}\frac{d}{dt}E^1_{(1-\beta,1),1}\left(\frac{-1}{K}(t-\tau)^{1-\beta},
   \frac{-\lambda_j}{K}(t-\tau)\right)\left(f(\tau),\varphi_{j}\right)d\tau+\left(f(\tau),\varphi_{j}\right).
\end{split}\end{equation}
By the complex contour integral representation of the bivariate Mittag-Leffler function in Eq. (\ref{est:M_LLL}), we get
\begin{equation}\begin{split}\label{eq:AAab}
E^1_{(1-\beta,1),1}\left(\frac{-1}{K}t^{1-\beta},
\frac{-\lambda_j}{K}t\right)
=&1-\frac{1}{K}t^{1-\beta} E_{(1-\beta,1),2-\beta}\left(\frac{-1}{K}t^{1-\beta},  \frac{-\lambda_j}{K}t\right)\\                                                  &-\frac{\lambda_j t}{K} E_{(1-\beta,1),2}\left(\frac{-1}{K}t^{1-\beta},  \frac{-\lambda_j }{K}t\right).
\end{split}\end{equation}
Combining Eq. (\ref{eq:AAa}) and Eq. (\ref{eq:AAab}), with Lemmata \ref{Lem:MandL} and \ref{Lemma:ML}, and Young's inequality, we have
\begin{align}\label{ieq:firstderiva}
&\left\|\partial_t u(t)\right\|^{2}_{L^{2}(\Omega)} \notag\\
=&\sum\limits_{j=1}^{\infty}\left(\frac{d}{dt}\int_{0}^{t}E^1_{(1-\beta,1),1}\left(\frac{-1}{K}(t-\tau)^{1-\beta},
   \frac{-\lambda_j}{K}(t-\tau)\right)\left(f(\tau),\varphi_{j}\right)d\tau\right)^{2}                    \notag\\
\leq &\sum\limits_{j=1}^{\infty}
      \left[\int_{0}^{t}\left(\frac{(t-\tau)^{-\beta}}{K}E^1_{(1-\beta,1),1-\beta}\left(\frac{-1}{K}(t-\tau)^{1-\beta},
           \frac{-\lambda_j}{K}(t-\tau)\right)\right.\right.                                                \notag\\                                                   &\left.\left.+C\sum\limits_{j=1}^{\infty}\int_{0}^{t}\left|\left(f(\tau),\varphi_{j}\right)\right|^2d\tau -\frac{\lambda_j}{K}E^1_{(1-\beta,1),1}\left(\frac{-1}{K}(t-\tau)^{1-\beta},
           \frac{-\lambda_j}{K}(t-\tau)\right)\right)\left|\left(f(\tau),\varphi_{j}\right)\right|d\tau\right]^2  \notag\\
\leq &C\sum\limits_{j=1}^{\infty}\left(\int_{0}^{t}\frac{(t-\tau)^{-\beta}}{K+\lambda_j(t-\tau)}d\tau\right)^2
      \int_{0}^{t}\left|\left(f(\tau),\varphi_{j}\right)\right|^2d\tau+C\sum\limits_{j=1}^{\infty}\int_{0}^{t}\left|\left(f(\tau),\varphi_{j}\right)\right|^2d\tau                                                  \notag\\
     &+C\sum\limits_{j=1}^{\infty}\left(\int_{0}^{t}\frac{(\lambda_j(t-\tau))^{\nu}(t-\tau)^{-\nu}}{K+\lambda_j(t-\tau)}d\tau\right)^2
      \int_{0}^{t}\left|\left(f(\tau),\varphi_{j}\right)\right|^2d\tau                                            \notag\\
\leq &C\left(t^{2-2\beta}+t^{2-2\nu}+1\right)\left\|f(t)\right\|_{L^{2}((0,T)\times\Omega)}
\leq C\left\|f(t)\right\|_{L^{2}((0,T)\times\Omega)},
\end{align}
where $\sup_{j\in\mathbb{N}}\frac{1}{K+\lambda_j(t-\tau)}\leq c$ and $\sup_{j\in\mathbb{N}}\frac{(\lambda_j(t-\tau))^{\nu}}{K+\lambda_j(t-\tau)}\leq c$  with $\nu\in(0,1)$ are used. Thus, $u(t)\in C^{1}([0,T],L^{2}(\Omega))$, and we have $\int_0^{T}\|\partial_tu(t)\|^{2}_{L^{2}(\Omega)}\leq c\|f(t)\|_{L^{2}((0,T)\times \Omega)}$.

With these analyses, we can deduce $u(t)\in C^{1}([0,T];L^2( \Omega))\cap L^2(0,T;H^{2}(\Omega)\cap H_0^{1}(\Omega))$ and $\partial_tu(t)\in L^{2}((0,T)\times\Omega)$.

In addition, with the help of
\begin{align*}
\mathfrak{L}_t\left\{ \partial_t^{\beta}\left[E^1_{(1-\beta,1),1}\left(\frac{-1}{K}t^{1-\beta},
     \frac{-\lambda_j}{K}t\right)\ast\left(f(\tau),\varphi_{j}\right)\right], z\right\}
= \frac{z^{-(1-\beta)}}{1+\frac{1}{K}z^{\beta-1}+\frac{\lambda_j}{K}z^{-1}}\left(\widehat{f}(z),\varphi_j\right),
\end{align*}
we deduce that
\begin{displaymath}\begin{split}
 &\partial_t^{\beta}\int_{0}^{t}E^1_{(1-\beta,1),1}\left(\frac{-1}{K}(t-\tau)^{1-\beta},
     \frac{-\lambda_j}{K}(t-\tau)\right)\left(f(\tau),\varphi_{j}\right)d\tau  \\
=&\int_{0}^{t}(t-\tau)^{-\beta}E^1_{(1-\beta,1),1-\beta}\left(\frac{-1}{K}(t-\tau)^{1-\beta},
     \frac{-\lambda_j}{K}(t-\tau)\right)\left(f(\tau),\varphi_{j}\right)d\tau.
\end{split}\end{displaymath}
Further, by Young's inequality, we obtain
\begin{align}\label{ieq:fracderiva}
\left\|\partial_t^{\beta}u(t)\right\|^{2}_{L^{2}(\Omega)}
=&\sum\limits_{j=1}^{\infty}\left(\int_{0}^{t}(t-\tau)^{-\beta}E^1_{(1-\beta,1),1-\beta}\left(\frac{-1}{K}(t-\tau)^{1-\beta},
   \frac{-\lambda_j}{K}(t-\tau)\right)\left(f(\tau),\varphi_{j}\right)d\tau\right)^{2}  \notag\\
\leq&C\sum\limits_{j=1}^{\infty}\left(\int_{0}^{t}\frac{(t-\tau)^{-\beta}}
      {\left(1+\frac{\lambda_j(t-\tau)}{K}\right)}\left|\left(f(\tau),\varphi_{j}\right)\right|d\tau\right)^{2}   \notag\\
\leq&Ct^{2-2\beta}\sum\limits_{j=1}^{\infty}\int_{0}^{T}\left|\left(f(\tau),\varphi_{j}\right)\right|^2dt
\leq C\left\|f(t)\right\|^2_{L^{2}((0,T)\times \Omega)}.
\end{align}
Now we have $\partial^\beta_t u(t)\in L^{2}((0,T)\times \Omega)$ and $\int_{0}^{T}\|\partial_t^{\beta}u(t)\|^{2}_{L^{2}(\Omega)}\leq C\|f(t)\|^2_{L^{2}((0,T)\times \Omega)}$. Thus, there is $\partial^\beta_tu(t)\in L^{2}((0,T)\times\Omega)$.

As for the well-posedness, according to Definition \ref{def:definition}, we only have to prove $\lim_{t\rightarrow 0}\|u(\cdot,t)\|_{\dot{H}^{-p}}=0$. In fact, by Eq. (\ref{inhomosolution}), there holds
\begin{align*}
\|u(\cdot,t)\|^2_{\dot{H}^{-p}}
=&\sum\limits_{j=1}^{\infty}\frac{1}{\lambda_j^{p}}\left(\int_{0}^{t}E^1_{(1-\beta,1),1}\left(\frac{-1}{K}(t-\tau)^{1-\beta},
   \frac{-\lambda_j}{K}(t-\tau)\right)\left(f(\tau),\varphi_{j}\right)d\tau\right)^{2}\\
\leq&\sum\limits_{j=1}^{\infty}\frac{1}{\lambda_j^{p}}\sup\limits_{0\leq\tau\leq T} \left|\left(f(\tau),\phi_j\right)\right|^{2}
\left|\int_{0}^{t}E^1_{(1-\beta,1),1}\left(\frac{-1}{K}(t-\tau)^{1-\beta},
   \frac{-\lambda_j}{K}(t-\tau)\right)d\tau\right|^2\\
\leq& C\|f(t)\|^2_{L^{\infty}(0,T;L^{2}(\Omega))} \sum\limits_{j=1}^{\infty}\frac{1}{\lambda_j^{p}}\left|\int_{0}^{t}E^1_{(1-\beta,1),1}\left(\frac{-1}{K}(t-\tau)^{1-\beta},
\frac{-\lambda_j}{K}(t-\tau)\right)d\tau\right|^2,
\end{align*}
where $\|f(t)\|_{L^\infty(0,T;L^2(\Omega))}:={\rm ess}\sup_{0\leq t\leq T}\|f(t)\|_{L^{2}(\Omega)}$. Since $\lambda_j\geq cj^{2/d}$ with $j\in \mathbb{N}_{+}$ (cf. \cite{Sakamoto11}), and $p>\frac{d}{2}$, then $\sum_{j=1}^\infty \frac{1}{\lambda_j^p}$ is bounded.
Also, with the boundness of the function $E^1_{(1-\beta,1),1}(\frac{-1}{K}(t-\tau)^{1-\beta},
   \frac{-\lambda_j}{K}(t-\tau))$, it follows that
\begin{displaymath}
\sum_{j=1}^{\infty}\frac{1}{\lambda_j^{p}}\left|\int_{0}^{t}E^1_{(1-\beta,1),1}(\frac{-1}{K}(t-\tau)^{1-\beta},
   \frac{-\lambda_j}{K}(t-\tau))d\tau\right|^2<\infty,
\end{displaymath}
and
$\lim\limits_{t\rightarrow 0}|\int_{0}^{t}E^1_{(1-\beta,1),1}(\frac{-1}{K}(t-\tau)^{1-\beta},
   \frac{-\lambda_j}{K}(t-\tau))d\tau|=0$ for each $j\in \mathbb{N}_{+}$.
With these, we know that $u(t)\in C([0,T];\dot{H}^{-p}(\Omega))$ with $p>d/2$. Also, by Lebesgue's dominant convergence theorem, there holds $\lim\limits_{t\rightarrow 0}\|u(\cdot,t)\|_{\dot{H}^{-p}}=0$. Therefore, $u(t)$ defined in Eq. (\ref{inhomosolution}) is indeed a weak solution to Problem (\ref{eq:problem}) with $u_0=0$. The uniqueness is the same as been proved in Theorem \ref{thm:regularityhomo}(I). Thus the proof is completed. \qed
\end{proof}

\subsubsection{Comparison of our results with standard diffusion equation and sub-diffusion equation}
Comparing with the classical diffusion equation (i.e., $K=0$ and $\beta=1$) and the sub-diffusion equation (i.e., $K=0$ and $0<\beta<1$), from the above analyses, the new model (i.e., $K>0$) discussed in this paper possesses the following interesting features:

\label{subsec:discuss}
\begin{description}
  \item[(i)] From the view of the physical background of the models, the classical diffusion equation is a macro description of Brownian motion; the sub-diffusion equation reflects that the transport rate of the micro particles is slower than that of Brownian motion; while the model (\ref{eq:problem}) with $K>0$  depicts a crossover from normal diffusion (as $t\rightarrow0$) to sub-diffusion (as $t\rightarrow\infty$) (see Appendix A for details). Therefore, the model (\ref{eq:problem}) discussed in this paper can improve and enrich the modeling delicacy in depicting the anomalous diffusion.

  \item[(ii)] In terms of the regularity of the solution for $f\equiv0$, if the initial value $u_0\in L^{2}(\Omega)$, then by Eq. (\ref{eq:nonsmoothest}), we have the same results as for the classical situation ($\|u(t)\|_{\dot{H}^{2}(\Omega)}$ $\leq C t^{-1}\|u_0\|_{L^{2}(\Omega)}$), but it is different from the case of sub-diffusion ($\|u(t)\|_{\dot{H}^{2}(\Omega)}\leq C t^{-\beta}\|u_0\|_{L^{2}(\Omega)}$, cf. \cite{Sakamoto11}), which has weak singularity w.r.t. the initial value.
If the initial value is smooth, i.e., $u_0\in H^{2}(\Omega)\cap H^0_1(\Omega)$, then from Theorem \ref{thm:regularityhomo} (III), the solution of the new model (\ref{eq:problem}) satisfies that $u(t)\in C^{1}([0,T],L^{2}(\Omega))\cap C([0,T];L^2(\Omega))\cap  C([0,T];H^2(\Omega)\cap H_0^1(\Omega))$, comparing with the corresponding results of the regularity of the sub-diffusion $u(t)\in C([0,T];L^2(\Omega))\cap C([0,T];H^2(\Omega)\cap H_0^1(\Omega))$ and $\partial_t^{\beta}u(t)\in C([0,T];L^2(\Omega))\cap C([0,T];H_0^1(\Omega))$ (cf. \cite{Sakamoto11,Jin18,Mu17}). Therefore, we deduce that the regularity of the solution to Problem (\ref{eq:problem}) with $K>0$ is indeed improved. With this, we conclude that the time regularity of the solution to TFDE is dominated by the highest-order operator in it.

  \item[(iii)] Theorem \ref{thm:regularityhomo} (I) shows the decay of the solution to Model (\ref{eq:problem}) is $t^{-1}$ as $t\rightarrow\infty$, which is similar to that of the classical diffusion problem. However, Corollaries 2.6 and 2.7 in \cite{Sakamoto11} show that the decay of solution to the sub-diffusion model is $t^{-\beta}$ as $t\rightarrow\infty$, which is  slower than these situations in classical diffusion equation as well as in new model.

  \item[(iv)] If $u_0\equiv0$, $f\in L^{2}((0,T)\times \Omega)$, and $\beta=1$, then the prior estimation for the sub-diffusion equation is (Theorem 2.2, Eq. (2.5) in \cite{Sakamoto11}):
      \begin{displaymath}
      \left\|u(t)\right\|_{L^{2}\left([0,T];\dot{H}^{2}(\Omega)\right)}+\|\partial^{\beta}_t u(t)\|_{L^{2}((0,T)\times \Omega)}\leq C\|f\|_{L^{2}((0,T)\times \Omega)},
      \end{displaymath}
      which is consistent with the estimation for the classical diffusion equation,
 i.e., \\ $\|u(t)\|_{L^{2}(0,T;\dot{H}^{2}(\Omega))}+\|\frac{\partial}{\partial t}u(t)\|_{L^{2}((0,T)\times \Omega)}\leq C\|f\|_{L^{2}((0,T)\times \Omega)}$ (see Page 382, Theorem 5 in \cite{Evans10}).
       This fact is also the same with our estimation in Eq. (\ref{sourceterm}) for $\beta=1$.
\end{description}

\section{Time semi-discretization and error estimation}
\label{sec:timediscrete}
In this section, we shall develop CIM for time semi-discretization. Error analysis will also be carried out.

\begin{figure}[t!]
\setlength{\abovecaptionskip}{0.1cm}  
\setlength{\belowcaptionskip}{-0.2cm} 
\centering
\begin{minipage}[c]{0.42\textwidth}
 \centering
 \centerline{\includegraphics[height=3.5cm,width=4.6cm]{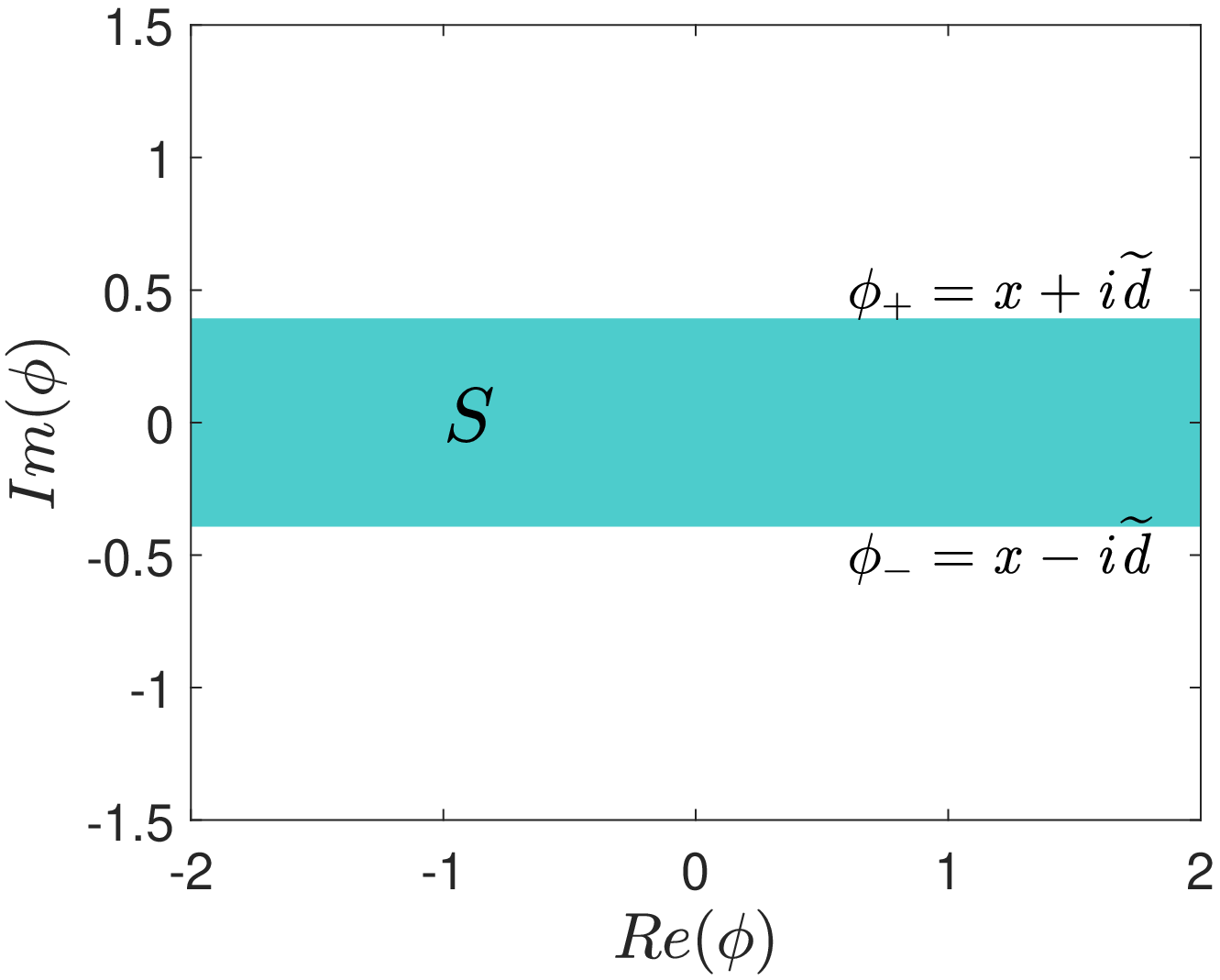}}
\end{minipage}
\begin{minipage}[c]{0.42\textwidth}
 \centering
 \centerline{\includegraphics[height=3.5cm,width=4.6cm]{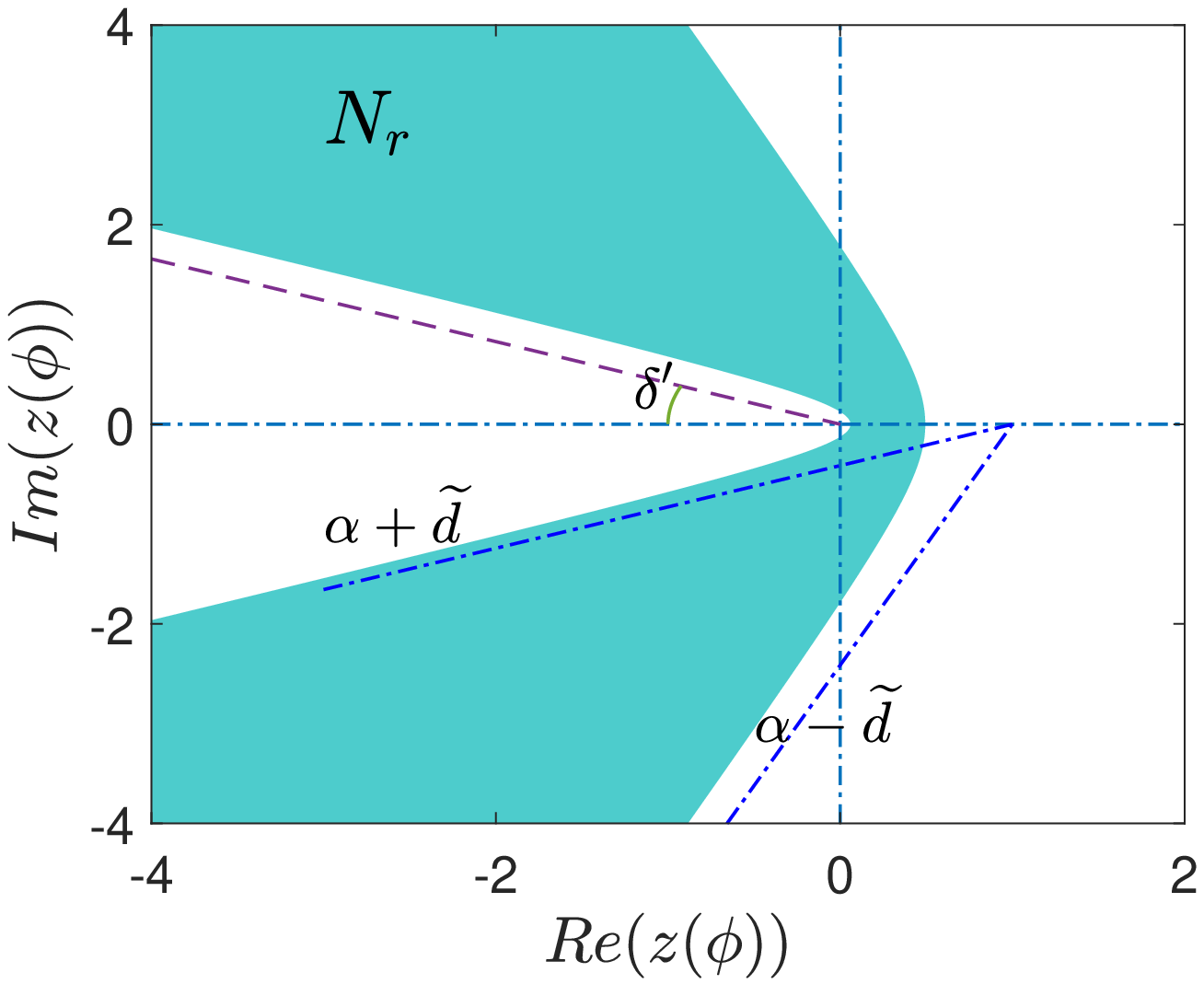}}
\end{minipage}
\caption{A diagram of the symmetrical open strip $S$ (\emph{left}) and the neighbourhood $N_r$ (\emph{right}).\emph{The parameters are taken as $\alpha=\pi/4$, $\mu=0.8$ and $\delta'=\pi/8$.}}
\label{Fig:example0}
\end{figure}

\subsection{CIM for Problem (\ref{eq:problem})}\label{subsection:3.1}
Following the basic ideas of CIM introduced in Sec. \ref{sec:Intro}, an appropriate integral contour must be selected firstly. In this paper, we choose the hyperbolic integral contour (cf. \cite{Fernandz04,Fernandz06}), which is parameterized as
\begin{equation}\label{eq:hyperboliccontour}
\Gamma:\ z(\phi) = \mu(1+\sin(i\phi-\alpha))~ ~ \quad {\rm with}~ \phi\in S,
\end{equation}
where $\mu>0$, $\alpha>0$ are parameters that need to be determined (see Eq. (\ref{para:Fernandz}));  $S$ is an open strip (see Figure \ref{Fig:example0} left), defined as $S:=\{\phi=x+iy\in\mathbb{C}, x\in\mathbb{R}, |y|< \tilde{d}\}$, where $\tilde{d}=\min\{\alpha, \pi/2-\alpha-\delta'\}$  with $0<\alpha-\tilde{d}<\alpha+\tilde{d}<\frac{\pi}{2}-\delta'$ and $\delta'\in(0,\pi/2)$ is the dip angle of asymptote of hyperbola (\ref{eq:hyperboliccontour}). The reason for designing such $S$ is given in the sequence.

\begin{remark}
    Indeed, the parabolic integral contour mentioned in \cite{Colbrook22a,Weideman07} is also a very efficient choice if all of the singularities of $\widehat{u}(z)$ are located on the negative real axi (see \cite{Weideman07} for details). Also, according to \cite{Weideman07}, large $N$ (the number of discretization points, which corresponding to the size of linear systems to be solved) would be used for the parabola than the hyperbola at the same required accuracy. Thus, in the spirit of efficiency, here we only consider the hyperbola one.
\end{remark}

Recalling the expression of $z(\phi)$ in Eq. (\ref{eq:hyperboliccontour}), the image of the horizontal line $\phi=x+iy\in S$ is
\begin{displaymath}
z(x+iy)=\mu\left(1-\sin(\alpha+y)\cosh(x)\right)+i\mu\cos(\alpha+y)\sinh(x),
\end{displaymath}
which can be expressed by the hyperbola
\begin{displaymath}
\left(\frac{\mu-u}{\sin(\alpha+y)}\right)^{2}-\left(\frac{v}{\cos(\alpha+y)}\right)^{2}=\mu^{2}\quad {\rm with}~ z=u+iv.
\end{displaymath}
With these, we can see that, for $y>0$, as $y$ increases from $0$ to $\pi/2-\alpha$,
the left branch of this hyperbola $z(\phi)$ will close and degenerate into the negative real axis; for, $y<0$, when $y$ decreases from $0$ to $-\alpha$, $z(\phi)$ will widen and become a vertical line. Since the integral contour we choose cannot degenerate into the real axis (see Remark \ref{thm:remark1}), therefore the included angle between the asymptotic line of the hyperbola and the real axis (denotes as $\delta'$), can not be equal to zero. In other words, the upper bound of $y$ is $\pi/2-\alpha-\delta'$ with $\delta' >0$, and the width of the strip $S$ satisfies $|y|< \tilde{d}$, with $\tilde{d}=\min\{\alpha, \pi/2-\alpha-\delta'\}$.

According to the definition of $z(\phi)$, it is clear that $z(\phi)$ is a conformal transformation, which maps the open strip $S$ into a neighbourhood in $z$-plan (see Fig. \ref{Fig:example0} (right)), and is denote as $N_r:=\{z(\phi)\in\mathbb{C}: \phi\in S\}$.  For a fixed $\delta'>0$, by Remark \ref{thm:remark1} and the definition of $z(\phi)$, there is $N_r\subseteq\Sigma_{\theta}$. Hence, the integrand $\widehat{u}(z)$ is analytic on $z(\phi)\in N_r$, and the integral (\ref{eq:integral}) can be written as
\begin{equation}\label{eq:integralsolve}
u(t)=I:=\int_{-\infty}^{+\infty} v(t,r,\phi)d\phi,
\end{equation}
where
\begin{equation}\label{exp:V1V2}
v(t,r,\phi)=\frac{1}{2\pi i}e^{z(\phi)t}\widehat{u}(z(\phi))z'(\phi),
\end{equation}
and $\widehat{u}(z)$ is defined in Eq. (\ref{eq:integrand}).

Assuming that the Laplace transform of the source term $f(t)$ exists (i.e., $|f(t)|<e^{\sigma_0 t}$), by Riemann-Schwarz reflection principle, there holds  $\widehat{f}(\overline{z})=\overline{\widehat{f}(z)}$. Since the integral contour $\Gamma$ defined in Eq. (\ref{eq:hyperboliccontour}) is symmetric w.r.t. the real axis, therefore  $\widehat{u}\left(\overline{z(\phi)}\right)=\overline{\widehat{u}(z(\phi))}$. After using the mid-point rule to approximate integral in Eq. (\ref{eq:integralsolve}) with uniform step-spacing $\tau$, we can obtain CIM for time semi-discretization of Problem (\ref{eq:problem}):
\begin{equation}\begin{split}\label{eq:computeu}
 u^{N}(t)= I_{\tau;N}:= \tau\sum\limits_{k = 1-N}^{N-1}v(t,r,\phi_{k})
 &= \frac{\tau}{2\pi i}\sum\limits_{k = 1-N}^{N-1}e^{z(\phi_k)t}~\widehat{u}(z(\phi_k))z'(\phi_k)   \\ &=\frac{\tau}{\pi}\mathrm{Im}\left\{\sum\limits_{k = 0}^{N-1}e^{z(\phi_k)t}~\widehat{u}(z(\phi_k))z'(\phi_k)\right\},
\end{split}\end{equation}
where $\phi_{k}=(k+1/2)\tau$, $k=0,1,\cdots, N-1$.

\subsection{Convergence analysis of CIM scheme (\ref{eq:computeu})}
In this subsection, we shall analyze the convergence of scheme (\ref{eq:computeu}). For this, we further assume that $f(t)\in L^{1}(\Omega)$, then $\widehat{f}(z)$ might be continued as a analytic function in $\Sigma_{\theta}$ with fixed $\theta\in(\pi/2,\pi)$. With this, the estimates can be expressed in terms of $\widehat{f}(z)$ rather than $f(t)$ (cf. \cite{McLean04}). We also define
\begin{equation}\label{norm_f}
\left\|\widehat{f}(z)\right\|_D:=\sup\limits_{z\in D}\left\|\widehat{f}(z)\right\|_{L^{2}(\Omega)},~ ~ D\subseteq\Sigma_{\theta}.
\end{equation}
In this paper, we take $D=N_r$.

By introducing the series $I_{\tau;N;\epsilon}=\tau\sum_{k=1-N}^{N-1}v(t,r,\phi_{k})(1+\epsilon_k)$ (which will be explained in the sequence) and $I_{\tau}:=\tau\sum_{k=-\infty}^{\infty}v(t,r,\phi_k)$, the error of CIM scheme (\ref{eq:computeu}) can be splitted into
\begin{align}\label{eq:EN}
E_{N}:=\left\|u(t)-u^{N}(t)\right\|_{L^{2}(\Omega)}
\leq &\left\|I-I_\tau\right\|_{L^{2}(\Omega)}+\left\|I_\tau-I_{\tau;N}\right\|_{L^{2}(\Omega)}
+\left\|I_{\tau;N}-I_{\tau;N;\epsilon}\right\|_{L^{2}(\Omega)} \notag\\
=&DE+TE+RE,
\end{align}
where $DE: =\|I-I_\tau\|_{L^{2}(\Omega)}$ is the discretization error,  $TE:=\|I_\tau-I_{\tau;N}\|_{L^{2}(\Omega)}$ is the truncation error and $RE:=\|I_{\tau;N}-I_{\tau;N;\epsilon}\|_{L^{2}(\Omega)}$ is usually defined as the round-off error.

Here, we should remark that the propagation of the round-off error when applying CIM has already been addressed and resolved for the hyperbolas by  M. Lopez-Fernandez et al. in \cite{Fernandz06}. As discussed in \cite{Fernandz06}, the reason why we introduce the series $I_{\tau;N;\epsilon}$ is mainly because that $\frac{1}{2\pi i}e^{z(\phi_k)t}\hat{u}(z(\phi_k)) z'(\phi_k)$ is approximated by $\frac{1}{2\pi i}e^{z_{k}t}~\hat{u}_kz'_{k}$, where $\{\hat{u}_k\}_{k=1-N}^{N-1}$ are provided by means of solving the linear system (\ref{sys:linear}) with prescribed accuracy $\epsilon$. That is, $\|\hat{u}_k-\hat{u}(z(\phi_k))\|_{L^{2}(\Omega)}<\epsilon_k$, where $\epsilon_k$s are the relative errors in the computed function values $\hat{u}(z(\phi_k))$ and $|\epsilon_k|<\epsilon$ for all $k=\pm1, \pm2,... ,\pm (N-1)$. The evaluations of the elementary functions involved in $v(t,r,\phi_{k})$ (such as $\exp(z_k t)$, $z'_k$ and $z_k$, etc.) turn out to be negligible compared to $\epsilon$ (see \cite{Fernandz06}). Therefore, $\|e^{z_{k}t}\hat{u}_kz'_{k}-e^{z(\phi_{k})t}\hat{u}(\phi_{k})z'(\phi_{k}))\|_{L^{2}(\Omega)}\leq \epsilon_k\|e^{z_{k}t}\hat{u}(z(\phi_k)) z'(\phi_k)\|_{L^{2}(\Omega)}$ and $I_{\tau;N;\epsilon}$ make sense.

To estimate (\ref{eq:EN}), we need to discuss the properties of the integrand $v(t,r,\phi_{k})$.

Firstly, for $z\in N_r\subseteq\Sigma_{\theta}$ with $|z|>\delta$ (where $\delta$ satisfies the condition in Proposition \ref{thm:prop1}),  by Proposition \ref{thm:prop1}, there exists a positive constant $C$, such that the integrand $\widehat{u}(z)$ defined in Eq. (\ref{eq:integrand}) satisfies
\begin{equation}\label{ieq:uhatinhomo}
\left\|\widehat{u}(z)\right\|_{L^{2}(\Omega)}
\leq C|z|^{-1}\left(\|u_0\|_{L^{2}(\Omega)}+\left\|\widehat{f}(z)\right\|_{N_r}\right).
\end{equation}
The following lemmas are also needed.

\begin{lemma}[cf. Lemma 1 in \cite{Fernandz04}]\label{et:lemma}
Set $L(x):=1+\big|\ln(1-e^{-x})\big|$, $x>0$. There hold
\begin{displaymath}
\int_{0}^{+\infty}e^{-\gamma\cosh(x)}dx\leq L(\gamma),~ ~ \gamma>0,
\end{displaymath}
and for $\sigma>0$,
\begin{displaymath}
\int_{\sigma}^{+\infty}e^{-\gamma\cosh(x)}dx\leq (1+L(\gamma))e^{-\gamma\cosh(\sigma)},~ ~\gamma>0.
\end{displaymath}
\end{lemma}
Obviously, $L(x)\rightarrow1$  as $x\rightarrow\infty$ and $L(x)\sim|\ln x|$ as $x\rightarrow0^{+}$.
Then the integrand $v(t,r,\phi_{k})$ has the properties shown in the following proposition.

\begin{lemma}\label{thm:prop2}
Let $v(t,r,\phi)$ be defined in Eq. (\ref{exp:V1V2}) with $z(\phi)\in N_r\subset\Sigma_{\theta}$ and $\phi=x\pm iy\in S$. Then, for $t>0$, $\left\|\widehat{f}(z)\right\|_{N_r}<\infty$ and $u_0\in L^{2}(\Omega)$, there exists a positive constant $C$, such that
\begin{equation}\label{ieq:setev}
\|v(t,r,\phi)\|_{L^{2}(\Omega)}
\leq C\varphi\left(\alpha,\tilde{d}\right)e^{\mu t}\left(\|u_0\|_{L^{2}(\Omega)}
+\left\|\widehat{f}(z)\right\|_{N_r}\right)e^{-\mu t\sin(\alpha)\cosh(x)},
\end{equation}
where $\varphi(\alpha,\tilde{d})=\sqrt{\frac{1+\sin(\alpha+\tilde{d})}{1-\sin(\alpha+\tilde{d})}}$.
\end{lemma}

Since one can obtain the above result by proceeding as in the proof of Theorem 2 in \cite{Fernandz04} and together with (\ref{ieq:uhatinhomo}), without much difficulty, we omit the proof.

Now we are ready to present the error estimates of the CIM.

\begin{theorem}\label{thm:fernandz}
Let $u(t)$ be the solution to Problem (\ref{eq:problem}), and $u^{N}(t)$ be given by Eq.  (\ref{eq:computeu}), respectively. Given $\Lambda>1$, $\alpha$, $\delta'$ with $0<\alpha,\delta'<\pi/2$, and the prescribed round-off error $\epsilon$. For $0<t_0\leq t\leq\Lambda t_0$ and $0<\varrho<1$, if  $\left\|\widehat{f}(z)\right\|_{N_r}<\infty$, 
then there uniformly holds
\begin{align*}
\left\|u(t)-u^{N}(t)\right\|_{L^{2}(\Omega)}
\leq& CL\left(\mu t_0\sin(\alpha-\tilde{d})\right)\varphi(\alpha,\tilde{d})\left(\|u_0\|_{L^{2}(\Omega)}+\left\|\hat{f}(z)\right\|_{N_r}\right)\\
&\times\left(\epsilon\cdot(\epsilon_N(\varrho))^{\varrho-1}+\frac{(\epsilon_N(\varrho))^{\varrho}}{1-\epsilon_N(\varrho)}\right),
\end{align*}
where $\epsilon_N(\varrho)
:=\exp(-\frac{2\pi \tilde{d}N}{a(\varrho)})$, $a(\varrho):=\cosh^{-1}\left(\frac{\Lambda}{(1-\varrho)\sin(\alpha-\widetilde{d})}\right)$, $0<\tilde{d}<\min\{\alpha,\pi/2-\alpha-\delta'\}$, and the parameter $\mu$ in Eq. (\ref{eq:hyperboliccontour}) is chosen as $\mu=\frac{2\pi \tilde{d}N(1-\varrho)}{t_0\Lambda a(\varrho)}$.
\end{theorem}
\begin{proof}

Since the integrand $v(t,r,\phi)$ with $\phi=x\pm iy\in S$ has the decay property described in Proposition \ref{thm:prop2},  by following the proof of Theorem 2 in \cite{Fernandz04}, we can get the uniform estimates of DE and TE on $t_0\leq t\leq t_1=\Lambda t_0$ with $t_0>0$ and $\Lambda=t_1/t_0>1$:
\begin{align}\label{eq:EN_11}
DE+TE\leq &C\varphi(\alpha,\tilde{d})L(\mu t_0\sin(\alpha-\tilde{d}))\left(\|u_0\|_{L^{2}(\Omega)}
+\left\|\widehat{f}(z)\right\|_{N_r}\right)        \notag\\
&\times e^{\mu t_0\Lambda}\left(\frac{1}{e^{2\pi \tilde{d}/\tau}-1}+\frac{1}{e^{\mu t_0\sin(\alpha-\tilde{d})\cosh(N\tau)}}\right).
\end{align}

Take
\begin{equation}\label{para:Fernandz}
\tau=\frac{a(\varrho)}{N},~~\textrm{and}~\mu=\frac{2\pi \tilde{d}N(1-\varrho)}{t_0\Lambda a(\varrho)}.
\end{equation}


Then by Theorem 1 in \cite{Fernandz06}, and with the help of Eq. (\ref{para:Fernandz}), we get
\begin{equation}\label{eq:ENN}
DE+TE\leq C\varphi(\alpha,\tilde{d})\cdot L(\mu t_0\sin(\alpha-\tilde{d}))\left(\|u_0\|_{L^{2}(\Omega)}
+\left\|\widehat{f}(z)\right\|_{N_r}\right)e^{\mu\Lambda t_0}\frac{2\epsilon_N(\varrho)}{1-\epsilon_N(\varrho)},
\end{equation}


As for $RE$, since
\begin{displaymath}
RE=\left\|I_{\tau;N}-I_{\tau;N;\epsilon}\right\|_{L^{2}(\Omega)}
\leq \tau\sum_{k=1-N}^{N-1}\left\|v(t,r,\phi_{k})\right\|_{L^{2}(\Omega)}|\epsilon_k|
\leq 2\epsilon\tau\sum\limits_{k=0}^{N-1}\left\|v(t,r,\phi_{k})\right\|_{L^{2}(\Omega)},
\end{displaymath}
then by Proposition \ref{thm:prop2} and Lemma \ref{et:lemma}, and taking $\tau$ and $\mu$ as in (\ref{para:Fernandz}), for $t_0\leq t\leq \Lambda t_0$, there holds
\begin{equation}\begin{split}\label{err:roundoff}
RE &\leq 2C\epsilon e^{\mu\Lambda t_0}\varphi(\alpha,\tilde{d})\left(\|u_0\|_{L^{2}(\Omega)}
 + \left\|\widehat{f}(z)\right\|_{N_r}\right)\tau\sum\limits_{k=0}^{N-1}e^{-\mu t_0\sin(\alpha-\tilde{d})\cosh (k+1/2)\tau}\\
&\leq CL(\mu t_0\sin(\alpha-\tilde{d}))\varphi(\alpha,\tilde{d})\left(\|u_0\|_{L^{2}(\Omega)}+\left\|\widehat{f}(z)\right\|_{N_r}\right)\epsilon e^{\mu\Lambda t_0}.
\end{split}\end{equation}
Note that (cf. \cite{Fernandz06})
\begin{equation}\label{remark:note}
e^{\mu\Lambda t_0}\epsilon_N(\varrho)=(\epsilon_N(\varrho))^{\varrho-1}\cdot \epsilon_N(\varrho)=(\epsilon_N(\varrho))^{\varrho}.
\end{equation}
Together above, we finally get the error estimation $E_N$ for CIM scheme (\ref{eq:computeu}).\qed

\end{proof}

\begin{corollary}\label{co:delta-solution}
Given $t_0>0$ and $\Lambda>1$. For $0< t\leq (\Lambda-1) t_0$ and $0<\varrho<1$, there exists a constant $C_1$, which depends on $\Omega,K,\beta,\alpha,\rho,N,\mu$, such that
\begin{align*}
\left\|\sum_{k=-(N-1)}^{N-1}e^{z(\phi_k)t}z'(\phi_k)\right\|_{L^{2}(\Omega)}
\leq& C_1 L\left(\mu t_0\sin(\alpha-\tilde{d})\right)\varphi(\alpha,\tilde{d})\cdot\left(\epsilon\cdot(\epsilon_N(\varrho))^{\varrho-1}+\frac{(\epsilon_N(\varrho))^{\varrho}}{1-\epsilon_N(\varrho)}\right),
\end{align*}
uniformly holds, where $\phi_{k}=(k+1/2)\tau$, $k=0,1,\cdots, N-1$, and the other parameters are given as in Theorem \ref{thm:fernandz}.
\end{corollary}

We put the corresponding proof in Appendix C.

For fixed $\alpha$, $\Lambda$ and $\delta'$,  the optimal free parameter $\varrho^*$ is determined by minimizing the convex function
\begin{displaymath}
\epsilon\cdot (\epsilon_N(\varrho))^{\varrho-1}+\frac{(\epsilon_N(\varrho))^{\varrho}}{1-\epsilon_N(\varrho)}.
\end{displaymath}
Then the optimization parameters $\mu^*$ and $\tau^*$ are determined by submitting $\varrho^*$ into Eq. ({\ref{para:Fernandz}}).
From these, we can obtain the convergence order of CIM as follows.
\begin{theorem}\label{EQ:convergence}
For given $\Lambda$, and $\alpha$, $\delta'$ with $0<\alpha, \delta'<\pi/2$ and $0<\tilde{d}<\min\{\alpha,\pi/2-\alpha-\delta'\}$, by choosing the parameters provided in Eq. (\ref{para:Fernandz}) as the optimal parameters of the hyperbolic contour $\Gamma$, the actual error behaves as
\begin{displaymath}
E_N=\mathcal{O}\left(\epsilon+e^{-CN}\right)~~{\rm with}~~C=\mathcal{O}\left(\frac{1}{\ln(\Lambda)}\right).
\end{displaymath}
\end{theorem}
\begin{proof}
Based on the previous error analysis, we have $E_N=\mathcal{O}(e^{\mu t_0\Lambda}\varepsilon+e^{\mu t_0\Lambda-2\pi \tilde{d}/\tau})$. For the optimal parameters determined in Eq. (\ref{para:Fernandz}), there holds $C:=\mu t_0\Lambda-2\pi \tilde{d}/\tau=-\frac{2\pi \tilde{d}N\varrho}{a(\varrho)}$ and $a(\varrho)>\cosh^{-1}(\Lambda)\sim \ln(\Lambda)$ with $\Lambda>1$. Therefore, $C=\mathcal{O}\left(\frac{1}{\ln(\Lambda)}\right)$. As for $e^{\mu t_0\Lambda}\varepsilon$, because $t_1=t_0\Lambda$ is a given (finite) number and $\mu>0$ is also finite for given $N$ and $\Lambda$, so we have $e^{\mu t_0\Lambda}\varepsilon=\mathcal{O}(\varepsilon)$. \qed
\end{proof}

Numerical results in Fig. \ref{Fig:example1} (right) and Fig. \ref{Fig:example2} (b) in the sequence (where $\epsilon=2.22\times 10^{-16}$ is chosen as the machine precision in IEEE arithmetic) demonstrate that the convergence order $C=\mathcal{O}\left(\frac{1}{\ln(\Lambda)}\right)$ in Theorem \ref{EQ:convergence} is indeed optimal. That is, the CIM has spectral accuracy w.r.t. $N$.

\begin{remark}
    We re-emphasize that, as mentioned in the end of Subsection \ref{subsection:3.1}, the error estimate results of the CIM scheme are independent of the smoothness of the solution (in the temporal direction), which is one of the typical advantages of the CIM. Theorem \ref{thm:fernandz} (or (\ref{exp:V1V2}), (\ref{eq:computeu}) and (\ref{eq:EN_11})) shows/show that high accuracy can be reached by using the CIM only if the real solution satisfies the conditions of Laplace transform, which is a very loose regularity requirement for almost all problems we can encounter.
\end{remark}

\begin{remark}
We note that: (a) there are some places where similar rates have been shown (with logarithmic dependence in $C$), e.g. \cite{Colbrook22a} and \cite{Colbrook22b}. More specifically, here the result we obtained in Theorem \ref{EQ:convergence} is consistent with that in \cite{Fernandz06}, i.e. $O(e^{-CN})$, while the convergence rate in \cite{Colbrook22a} and \cite{Colbrook22b} is $O(e^{-cN/\log(N)})$. (b) as mentioned in \cite{Weideman10}, the round-off error is an important factor in the stability issue of the CIM. In our error estimates, we have taken this into account. So the corresponding CIM scheme is stable and does not have the situation of error inverse growth, which will be also visually demonstrated in the section of numerical experiments.
\end{remark}

\section{Spatial semi-discretization and error estimation}
\label{sec:Galerkin}
In this section, we develop the spatial semi-discrete scheme by using the standard Galerkin finite element method, and give error estimates for the semi-discrete scheme of space with both smooth and non-smooth initial values, respectively.
\subsection{Semi-discrete Galerkin scheme}
Let $\mathcal{T}_{h}=\{_{\mathcal{T}_{j}}\}_{j=1}^{M}$ belong to a family of quasi-uniform triangulations of $\Omega$ with $h = \max_{0\leq j\leq M}{\rm diam}(\tau_{j})$, where the boundary triangles are allowed to have one curved edge along $\partial\Omega$, and let $S_{h}$ denote the piecewise continuous functions on the closure $\overline{\Omega}$ of $\Omega$ which reduce to polynomials of degree $\leq r-1$ on each triangle $\tau_j$ and vanish outside $\Omega_{h}$, namely,
\begin{displaymath}
 S_{h}=\left\{\chi\in H_{0}^{1}(\Omega)\cap \mathrm{C}\left(\overline{\Omega}\right); \chi\mid_{\mathcal{T}_{j}}\in\Pi_{r-1}\right\},
\end{displaymath}
where $\Pi_{s}$ denotes the set of polynomials of degree at most $s$.
Then the $L^{2}$-orthogonal projection $P_{h}: L^{2}(\Omega)\rightarrow S_{h}$ can be defined as
\begin{displaymath}
\left(P_{h}u,v_{h}\right)=\left(u,v_{h}\right)\quad \forall~v_{h}\in S_{h},~ u\in L^{2}(\Omega),
\end{displaymath}
and the standard  Ritz  projection $R_{h}:H_{0}^{1}(\Omega)\rightarrow S_{h}$ is defined as
\begin{displaymath}
\left(\nabla R_{h}u,\nabla v_{h}\right)=\left(\nabla u,\nabla v_{h}\right)\quad \forall~v_{h}\in S_{h},~u \in H_{0}^{1}(\Omega).
\end{displaymath}
For $q=1$, $2$, it is well known that the  projections $P_{h}$ and  $R_{h}$ meet the following estimates (see \cite{Thomee06}, etc.), respectively,
\begin{equation}\label{eq:PQ}
\begin{split}
\left\|P_{h}v-v\right\|_{L^{2}(\Omega)}+h\left\|\nabla(P_{h}v-v)\right\|_{L^{2}(\Omega)} & \leq Ch^{q}\|v\|_{\dot{H}^{q}(\Omega)}~ ~ {\rm for~ all}~v\in \dot{H}^{q}(\Omega),\\
\left\|R_{h}v-v\right\|_{L^{2}(\Omega)}+h\left\|\nabla(P_{h}v-v)\right\|_{L^{2}(\Omega)} & \leq Ch^{q}\|v\|_{\dot{H}^{q}(\Omega)}~ ~ {\rm for~ all}~v\in \dot{H}^{q}(\Omega).
\end{split}
\end{equation}
Besides, $P_h$ satisfies the stability estimate $\|P_{h}v\|_{L^{2}(\Omega)}\leq C\|v\|_{L^{2}(\Omega)}$ for $v\in L^{2}(\Omega)$.

Based on above, the spatially semi-discrete Galerkin FEM scheme for Problem (\ref{eq:problem}) can be formally described as: for $t>0$, find $u_h(t)\in S_h$, such that
\begin{equation}\label{eq:finitele}
\left(K\frac{\partial u_h(t)}{\partial t}, v_h\right)+\left(\partial_{t}^{\beta}u_h(t), v_h\right)+(Au_h(t),v_h)=(f,v_h)\quad \forall~v_h\in S_h
\end{equation}
with the initial value $u_h(0)=u_{0,h}$. The choice of $u_{0,h}$ depends on the smoothness of the initial data $u_0$, i.e., we take $u_{0,h}=R_hu_0$ for $u_0\in H_{0}^{1}(\Omega)\cap H^{2}(\Omega)$, and take $u_{0,h}=P_h u_0$ when $u_0\in L^{2}(\Omega)$.
In addition, by introducing the finite element version of the operator $A$ as $A_{h}: S_{h}\rightarrow S_{h}$ with
\begin{displaymath}
\left(A_{h}u_{h}(t),v_{h}\right)=\left(\nabla u_{h}(t), \nabla v_{h}\right)\quad \forall\ u_{h}(t),~ v_{h}\in S_{h},
\end{displaymath}
then Eq. (\ref{eq:finitele}) turns to:
\begin{equation}\label{eq:semispace}
K\frac{\partial u_h(t)}{\partial t}+\partial_{t}^{\beta}u_h(t)+A_hu_h=f_h\quad \forall~ t>0
\end{equation}
in the weak sense, with $f_h(t)=P_hf(t)$. After taking the Laplace transform on both sides of Eq. (\ref{eq:semispace}), there holds
\begin{equation}\label{eq:Lpsemispace}
\left(\left(Kz+z^{\beta}\right)I+A_h\right)\widehat{u}_h(z)=\left(K+z^{\beta-1}\right)u_{0,h}+\widehat{f}_{h}(z).
\end{equation}
Transforming the above equation by the inverse Laplace, and denote $W_h(z):=\frac{\eta(z)}{z}(\eta(z)I+A_h)^{-1}$,  we have
\begin{equation}\label{eq:spacesolut}
u_h(t)=\frac{1}{2\pi i}\int_{\Gamma_{\theta,\delta}}e^{zt}\left(W_h(z)u_{0,h}+z\eta^{-1}(z)W_h(z)\widehat{f}_{h}(z)\right)dz,
\end{equation}
where $W_h(z)$ satisfies the following estimates by analogy with Proposition \ref{thm:prop1} (II).
\begin{corollary}\label{lemma:w_h}
Let $\delta$ satisfy the condition in the statement of Proposition \ref{thm:prop1}. For all $z\in\Sigma_{\theta,\delta}$, there hold
\begin{equation}
\left\|W_h(z)\right\|_{L^{2}(\Omega)}\leq C|z|^{-1}, \quad \left\|A_hW_h(z)\right\|_{L^{2}(\Omega)}\leq C.
\end{equation}
\end{corollary}

\subsection{The spatial error estimation}
In this subsection, we shall perform the error estimation of scheme (\ref{eq:semispace}). For the homogeneous case,  the spatial error estimation in $L^{2}$-norm with both smooth and non-smooth initial values are shown as follows.
\begin{theorem}\label{th:spaceerrorhomo1}
Let $u(t)$ be the solution to Problem (\ref{eq:problem})  and $u_{h}(t)$ be given in Eq. (\ref{eq:spacesolut}), respectively. If $f\equiv0$, then for $0<t\leq T$, the following error estimations hold:
\begin{description}
  \item[(I)] Assume $u_0\in L^{2}(\Omega)$. Take $u_{0,h}=P_hu_{0}$. There holds
\begin{displaymath}
\left\|u_h(t)-u(t)\right\|_{L^{2}(\Omega)}\leq C t^{-1}h^{2}\|u_0\|_{L^{2}(\Omega)}.
\end{displaymath}
  \item[(II)] If $u_0\in H_{0}^{1}(\Omega)\cap H^{2}(\Omega)$, and $u_{0,h}=R_hu_{0}$, then there exists a positive constant $C$ such that
\begin{displaymath}
\left\|u_h(t)-u(t)\right\|_{L^{2}(\Omega)}\leq Ch^{2}\|u_0\|_{\dot{H}^{2}(\Omega)}.
\end{displaymath}
\end{description}
\end{theorem}

\begin{proof}

Let $f\equiv0$. For $t>0$ and $u_{0,h}=P_hu_{0}$, the error e(t):=$u(t)-u_h(t)$ is given by
\begin{align*}
e(t)
&=\frac{1}{2\pi i}\int_{\Gamma_{\theta,\delta}}e^{zt}(W_h(z)-W(z))u_0dz \\
&=\frac{1}{2\pi i}\int_{\Gamma_{\theta,\delta}}e^{zt}\frac{\eta(z)}{z}\left((\eta(z)I+A_h)^{-1}P_h-(\eta(z)I+A)^{-1}\right)u_0dz.
\end{align*}
 Denote $H:=(\eta(z)I+A_h)^{-1}P_h-(\eta(z)I+A)^{-1}$. It can be rewritten as
\begin{align*}
H =&P_h\left((\eta(z)I+A_h)^{-1}-(\eta(z)I+A)^{-1}\right)P_h \\
   &-(I-P_h)(\eta(z)I+A)^{-1}P_h-(\eta(z)I+A)^{-1}(I-P_h).
\end{align*}
Note that $A\left(\eta(z)I+A\right)^{-1}=I-\eta(z)(\eta(z)I+A)^{-1}$ is uniformly bounded for $z\in\Sigma_{\theta,\delta}$, which also holds with $A$ replaced by $A_h$. With this, according to the proof of Theorem 2.1 in \cite{Lubich96}, the operator $H$ is uniformly bounded for $z\in\Sigma_{\theta,\delta}$, and further we have $\|Hu_0\|_{L^{2}(\Omega)}\leq Ch^2\|u_0\|_{L^{2}(\Omega)}$. Therefore, we obtain
\begin{displaymath}
\left\|e(t)\right\|_{L^{2}(\Omega)}
\leq  Ch^2\int_{\Gamma_{\theta,\delta}}e^{|z|t\cos(\theta)}|dz|\|u_0\|_{L^{2}(\Omega)}\leq Ct^{-1}h^2\|u_0\|_{L^{2}(\Omega)}.
\end{displaymath}
Thus, the proof of Theorem \ref{th:spaceerrorhomo1} (I) is completed.

Let $u_0\in H_{0}^{1}(\Omega)\cap H^{2}(\Omega)$, $u_{0,h}=R_hu_{0}$ and $f\equiv0$. Similar to the analysis above, the  error $\tilde{e}(t):=u(t)-u_h(t)$ is given by
\begin{align*}
\tilde{e}(t)
&=\frac{1}{2\pi i}\int_{\Gamma_{\theta,\delta}}e^{zt}(W_h(z)R_h-W(z))u_0dz \\
&=\frac{1}{2\pi i}\int_{\Gamma_{\theta,\delta}}e^{zt}\frac{\eta(z)}{z}\left((\eta(z)I+A_h)^{-1}R_h-(\eta(z)I+A)^{-1}\right)u_0dz.
\end{align*}
Since $\eta(z)(\eta(z)I+A)^{-1}=I-(\eta(z)I+A)^{-1}A$, which also holds with $A$ replaced by $A_h$, then we have
\begin{align*}
\tilde{e}(t)=&\frac{1}{2\pi i}\int_{\Gamma_{\theta,\delta}}e^{zt}z^{-1}(R_h-I)u_0dz  \\
&-\frac{1}{2\pi i}\int_{\Gamma_{\theta,\delta}}e^{zt}z^{-1}\left((\eta(z)I+A_h)^{-1}A_hR_h-(\eta(z)I+A)^{-1}A\right)u_0dz.
\end{align*}
Since $A_hR_h=P_hA$, then we have $\left((\eta(z)I+A_h)^{-1}A_hR_h-(\eta(z)I+A)^{-1}A\right)u_0=HAu_0$. Further, by the boundness of $H$ and the estimate in Eq. (\ref{eq:PQ}), there holds
\begin{displaymath}
\left\|\tilde{e}(t)\right\|_{L^{2}(\Omega)}
\leq Ch^2\int_{\Gamma_{\theta,\delta}}e^{|z|t\cos(\theta)}|z|^{-1}|dz|\|Au_0\|_{L^{2}(\Omega)}\leq Ch^2\|u_0\|_{\dot{H}^{2}(\Omega)}.
\end{displaymath}
The proof is finished. \qed
\end{proof}

Next, we turn to the inhomogeneous case with vanishing initial data.
\begin{theorem}\label{th:spaceerrorinhomo}
Let $u(t)$ be the solution to Problem (\ref{eq:problem}) and $u_{h}(t)$ be given in Eq. (\ref{eq:spacesolut}), respectively. If $u_0\equiv0$, and $\left\|\hat{f}(z)\right\|_{N_r}<\infty$, then for $0<t\leq T$, there holds
\begin{displaymath}
\|u(t)-u_h(t)\|_{L^{2}(\Omega)}\leq C h^{2}t^{-1}\left\|\widehat{f}(z)\right\|_{N_r}.
\end{displaymath}
\end{theorem}

\begin{proof}
Since $u_0\equiv0$, the error $\bar{e}(t):=u_h(t)-u(t)$ is
\begin{align*}
\bar{e}(t)
&=\frac{1}{2\pi i}\int_{\Gamma_{\theta,\delta}}e^{zt}z\eta(z)^{-1}(W_h(z)P_h-W(z))\hat{f}(z)dz\\
&=\frac{1}{2\pi i}\int_{\Gamma_{\theta,\delta}}e^{zt}\left((\eta(z)I+A_h)^{-1}P_h-(\eta(z)I+A)^{-1}\right)\hat{f}(z)dz.
\end{align*}
According to the proof of Theorem \ref{th:spaceerrorhomo1} (I), it follows that
\begin{displaymath}
\|\bar{e}(t)\|_{L^{2}(\Omega)}\leq Ch^2\int_{\Gamma_{\theta,\delta}}e^{t|z|\cos(\theta)}\left\|\hat{f}(z)\right\|_{N_r}|dz|\leq C h^2 t^{-1}\left\|\widehat{f}(z)\right\|_{N_r}.
\end{displaymath} \qed
\end{proof}

\section{Fully discrete scheme and convergence analysis}
\label{sec:fulldiscret}
The fully discrete scheme of Problem (\ref{eq:problem}) can be obtained by applying CIM scheme to Eq. (\ref{eq:spacesolut}), which is given by
\begin{equation}\label{eq:fullcomputeu}
 u^{N}_{h}(t)= I_{\tau;h;N}:= \tau\sum\limits_{k = 1-N}^{N-1}v_h(t,r,\phi_{k})=\frac{\tau}{\pi}\mathrm{Im}\left\{\sum\limits_{k = 0}^{N-1}e^{z_{k}t}\widehat{u}_h(z_{k})z'_{k}\right\},
\end{equation}
where $\widehat{u}_h(z_k)$ is obtained by solving Eq. (\ref{eq:Lpsemispace}) with $z=z_k$, $k=0,1,2,..., N-1$. We can see that in order to get the numerical solution $u^{N}_{h}(t)$ for a range values of time $t$, one has to firstly solve $N$ discrete elliptic systems: $\forall v_{h}\in S_{h}$, $k=0,1,2,..., N-1$,
\begin{equation}\label{sys:linear}
\left(\left(Kz_k+z^{\beta}_k\right)\widehat{u}_h(z_k),v_h\right)+\left(\nabla \widehat{u}_{h}(z_k), \nabla v_{h}\right)=\left(\left(K+z^{\beta-1}_k\right)u_{0}+\widehat{f}(z_k),v_h\right).
\end{equation}
It is clear that the main computation of the fully discrete scheme (\ref{eq:fullcomputeu}) is to solve the  linear systems (\ref{sys:linear}). Moreover, all of the systems in (\ref{sys:linear}) are tri-diagonal linear algebraic equations with the coefficient matrix being strictly diagonally dominant. Hence, they have unique solutions.

In addition, by combining Theorem \ref{thm:fernandz} and Theorem \ref{th:spaceerrorhomo1}, we can directly get the following error estimation for the fully discrete scheme (\ref{eq:fullcomputeu}) with $f\equiv0$.
\begin{theorem}\label{thm:homog}
Let $u(t)$ be the solution to Problem (\ref{eq:problem}) and $u^{N}_{h}(t)$ be given by Eq. (\ref{eq:fullcomputeu}), respectively. For $t_0\leq t\leq\Lambda t_0$, $0<\varrho<1$ and $e_h^{N}:=u(t)-u^{N}_{h}(t)$, if $u_0\in L^{2}(\Omega)$,   there holds
\begin{align*}
\left\|e_h^{N}\right\|_{L^{2}(\Omega)}
\leq C\left(t^{-1}h^{2}
+ L(\mu t_0\sin(\alpha-\tilde{d}))\varphi(\alpha,\tilde{d})\left(\epsilon(\epsilon_N(\varrho))^{\varrho-1}
     +\frac{(\epsilon_N(\varrho))^{\varrho}}{1-\epsilon_N(\varrho)}\right)\right)\|u_0\|_{L^{2}(\Omega)}.
\end{align*}
If $u_0\in H^{2}(\Omega)\cap H^{1}_{0}(\Omega)$, we have
\begin{align*}
\left\|e_h^{N}\right\|_{L^{2}(\Omega)}
\leq C\left(h^{2}
+L(\mu t_0\sin(\alpha-\tilde{d}))\varphi(\alpha,\tilde{d})\left(\epsilon(\epsilon_N(\varrho))^{\varrho-1}
     +\frac{(\epsilon_N(\varrho))^{\varrho}}{1-\epsilon_N(\varrho)}\right)\right)\|u_0\|_{\dot{H}^{2}(\Omega)}.
\end{align*}
The parameters satisfy $0<\alpha,\delta'<\pi/2$ and $0<\tilde{d}<\min\{\alpha, \pi/2-\alpha-\delta'\}$.
\end{theorem}

Also, by combining Theorem \ref{thm:fernandz} and Theorem \ref{th:spaceerrorinhomo}, we can get the error estimation of the fully discrete scheme for the inhomogeneous problem with vanishing initial data, which is shown in the following theorem.
\begin{theorem}\label{Tm:fulldiscrete2}
Let $u(t)$ be the solution to Problem (\ref{eq:problem}) and $u^{N}_{h}(t)$ be given in Eq. (\ref{eq:fullcomputeu}), respectively. Given $\Lambda$, and $\alpha$, $\delta'$ with $0<\alpha, \delta'<\pi/2$, and $0<\tilde{d}<\min\{\alpha,\pi/2-\alpha-\delta'\}$, then for $t_0<t<\Lambda t_0$, $0<\varrho<1$ and $0<\gamma<1$ and $e_h^{N}:=u(t)-u^{N}_{h}(t)$, if $\left\|\widehat{f}(z)\right\|_{N_r}<\infty$, there holds
\begin{align*}
\left\|e_h^{N}\right\|_{L^{2}(\Omega)}
\leq C\left(t^{-1}h^{2}+L(\mu t_0\sin(\alpha-\tilde{d}))\varphi(\alpha,\tilde{d})\left(\epsilon (\epsilon_N(\varrho))^{\varrho-1}
     +\frac{(\epsilon_N(\varrho))^{\varrho}}{1-\epsilon_N(\varrho)}\right)
 \right)\left\|\widehat{f}(z)\right\|_{N_r}.
\end{align*}
\end{theorem}

\section{Implementation and acceleration of the algorithm}
\label{sec:acceler}
To compute one single linear elliptic system in (\ref{sys:linear}), our implementation process is as follows: for 1-D problem, the linear system is solved by Tomas algorithm, which has a computational complexity of $\mathcal{O}(M)$ (So that after combing with CIM algorithm for time discretization, the total computation amount for the solution to Problem (\ref{eq:problem}) is $\mathcal{O}(MN)$); for 2-D case, we directly employ the command $ mldivide(\cdot,\cdot)$ in MATLAB to solve the systems (\ref{sys:linear}).

In the sequence, we target on developing an acceleration algorithm for the fully discrete scheme (\ref{eq:fullcomputeu}).

It is clear that, given a contour line $z(\phi)$, to obtain the numerical solution from (\ref{eq:fullcomputeu}), the main calculation comes from getting $\{\widehat{u}_h(z_k)\}_{k=0}^{N-1}$ by solving $N$ elliptic equations in (\ref{sys:linear}). The computational cost will naturally by reduced a lot if we cut down the number of the elliptic systems.
Then the question is, how to keep the accuracy of the algorithm at the same time? The answer is that, since $\{\widehat{u}_h(z_k)\}_{k=0}^{N-1}$ used in (\ref{eq:fullcomputeu}) are mutually independent, they can be approximately interpolated by the reduced numbers of solutions of the elliptic systems in (\ref{sys:linear}) on some suitable points chosen from the contour line $z(\phi)$.
To avoid Runge phenomenon, we choose $n$ Chebyshev points as the corresponding interpolation nodes on the straight line $\phi$.

In numerical implementation, the specific contour $z(\phi)$ we select is the one
whose original image be just the real axis $\phi^0:=x+iy, x\in\mathbb{R}, y=0$.
As a matter of fact, by symmetry, we only need to discuss the details on half of $z(\phi^0)$ or $\phi^0$ (see Section \ref{sec:timediscrete} or (\ref{eq:fullcomputeu})).
\begin{figure}[t!]
\setlength{\abovecaptionskip}{0cm}  
\setlength{\belowcaptionskip}{-0.2cm} 
  \centering
  \includegraphics[height=5.5cm,width=10.0cm]{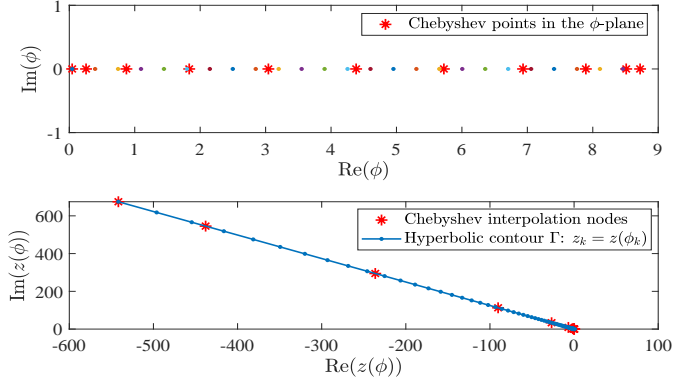}\\
\caption{ The Chebyshev interpolation nodes with $n=10$ and $N=100$. \emph{(top) the Chebyshev nodes (red-star) and the uniform partition points (color points, and four nodes are omitted between each two of them for clarity) in $\phi$-plane truncated at $N\tau^{*}$, where $\tau^{*}$ is determined by Eq. (\ref{para:Fernandz})); (bottom) uniform partition points and Chebyshev interpolation nodes in z-plane mapped by the hyperbolic contour $\Gamma$: $z(\phi)$) of (\ref{eq:hyperboliccontour}).}}
  \label{Fig:chebshev}
\end{figure}

Denote $n$ the Chebyshev points on the standard interval $[-1,1]$ as $\zeta_j=\cos(\frac{j\pi}{n})$, $j=0,1,2,...,n$, and $n< N$. So the corresponding Chebyshev interpolation nodes used are $z(x_j):=z(\frac{a+b}{2}+\frac{b-a}{2}\zeta_j))$, $j=0,1,2,...,n$ with $a=\frac{\tau}{2}$, $b=\left((N-1)+\frac{1}{2}\right)\tau$ (see Fig. \ref{Fig:chebshev}). We also denote the function values obtained by solving Eq. (\ref{sys:linear}) on the interpolation nodes ${z}(x_j)$ as $\widehat{u}_h({z}_{x_j})$, $j=0,1,2,...,n$.

With this, we select the barycentric Lagrange interpolation to approximate all $\{\widehat{u}_h(z_k)\}_{k=0}^{N-1}$. (The barycentric interpolation is a variant of Lagrange polynomial interpolation with $O(n)$ flops, which is fast and stable (see \cite{Trefethen04} for details).) That is
\begin{equation}\label{eq:interpolation}
\widehat{u}_h(z_k)\approx\widehat{u}^{n}_{I,h}(z_k):=\frac{\sum_{j=0}^{n}\frac{\omega_j}{z_k-{z}_{x_j}}\widehat{u}_h({z}_{x_j})}
{\sum_{j=0}^{n}\frac{\omega_j}{z_k-{z}_{x_j}}},~ ~
\omega_j=\frac{1}{\prod_{k\neq {m_j}}({z}_{x_j}-z_k)}, ~ j=0,1,...,n,
\end{equation}
where the barycentric weights $\omega_j$ are given as
\begin{displaymath}
\omega_j=
(-1)^{-j}\delta_j~ ~ {\rm with}~ ~ \delta_j=
\left\{
\begin{aligned}
&1/2,&& ~ ~ j=0~ {\rm or}~ j=n,\\
&1,&& ~ ~ {\rm otherwise}.
\end{aligned}\right.
\end{displaymath}
Until now, the fully discrete scheme (\ref{eq:fullcomputeu}) can be approximately replaced by
\begin{equation}\label{eq:fullcomputeu12}
 u^{N,n}_{I,h}(t):= \tau\sum\limits_{k = 1-N}^{N-1}v^{n}_{I,h}(t,\phi_{k})=\frac{\tau}{\pi}\mathrm{Im}\left\{\sum\limits_{k = 0}^{N-1}e^{z_{k}t}\widehat{u}^{n}_{I,h}(z_{k})z'_{k}\right\}.
\end{equation}

Combined with the spatial Galerkin FEM, the above process is shown in Algorithm \ref{CIM-FEM} with $\mathcal{O}(n(N+M))$ flops, which we name as ``CIM-Int-FEM''. In comparison, we also name the corresponding algorithm without acceleration as ``CIM-FEM''.

Denote $\phi_N^{+}:=[\frac{\tau}{2},\left((N-1)+\frac{1}{2}\right)\tau]$. We directly employ the result given in \cite{Trefethen04} on the standard interval $[-1,1]$ to explain the efficiency of our acceleration algorithm.   Firstly, we transform the interpolation from the complex plane that $z(\phi)$ lies in to the standard interval $[-1,1]$. Denote
$\widetilde{u}_h(\zeta):=\widehat{u}_h^N(z)$, $\widetilde{u}_{I,h}(\zeta):=\widehat{u}^{N,n}_{I,h}(z)$, $\zeta\in[-1,1]$, where $\widehat{u}_h^N(z)$ and $\widehat{u}^{N,n}_{I,h}(z)$ are the truncation of functions $\widehat{u}_h(z)$ and $\widehat{u}^{n}_{I,h}(z)$ on $z(\phi_N^{+})$, respectively.
Then by the discussion in Section 6 of \cite{Trefethen04}, the remainder term ${R}_n(z):=\widehat{u}^{N}_h(z)-\widehat{u}^{N,n}_{I,h}(z)$ satisfies the following estimation:
\begin{equation}\label{Eq:K}
\max_{z\in \phi_N^{+}} \left|R_n(z)\right|
=\max_{z\in \phi_N^{+}}\left|\widehat{u}_h^N(z)-\widehat{u}^{N,n}_{I,h}(z)\right|
=\max_{\zeta\in[-1,1]}\left|\widetilde{u}_h(\zeta)-\widetilde{u}_{I,h}(\zeta)\right|
\leq C \tilde{K}^{-n},
\end{equation}
where $\tilde{K}>1$ is a constant, determined by the continued analyticity of $\widetilde{u}_h(\zeta)$ in the complex plane (see Page 174 in \cite{Fornberg96} for more details).
To be precise, if $\widetilde{u}_h(\zeta)$ can be analytically continued to a function $\widetilde{u}_h$ in a complex region around $[-1,1]$ that includes an ellipse with foci $\pm1$ and axis lengths $2L$ and $2l$, then we may take $\tilde{K} = L + l$ (see Chapter 5 in \cite{Trefethen00} or Section 6 in \cite{Trefethen04}). In the sequence, we shall demonstrate the exact values of $L$ and $l$ taken in our situation.

\begin{figure}[t!]
\setlength{\abovecaptionskip}{0cm}  
\setlength{\belowcaptionskip}{-0.2cm} 
\centering
\begin{minipage}[c]{0.39\textwidth}
 \centerline{\includegraphics[height=5.5cm,width=10.6cm]{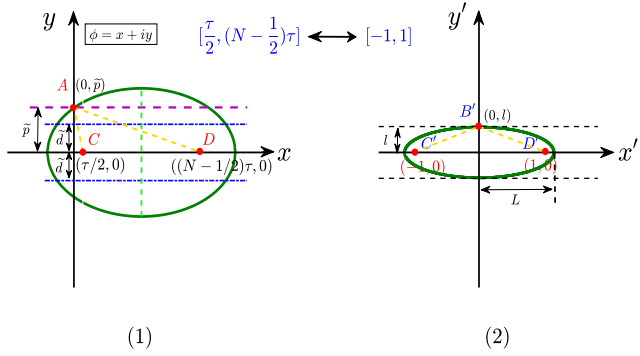}}
\end{minipage}
\caption{Schematic diagram of the maximal elliptic analytic region and its degradation to an elliptic region with focus $\pm1$, where $A:(0,\frac{\pi}{2}-\alpha-\varepsilon)$,
$C:(\frac{\tau}{2},0)$,
$D:\left(\left(N-\frac{1}{2}\right)\tau,0\right)$, $B':\left(0,l\right)$, $C':(-1,0)$, $D':(1,0)$.}
\label{Fig:example222}
\end{figure}

In fact, we only need to fix the biggest corresponding ellipse in the complex plane that $\phi_N^+$ lies in. As we have discussed in Subsection \ref{subsection:3.1}, $\widehat{u}_h(z(\phi_N^{+}))$ can be analytically continued to $\widehat{u}_h(z(\phi))$ in the $\phi$-complex plane until $z(\phi)=0$ or $\phi=i\cdot\left(\frac{\pi}{2}-\alpha+2m\pi\right), m\in \mathbb{Z}$. So the biggest ellipse, on and inside which $\widehat{u}$ is analytic, with foci $C=(\frac{\tau}{2},0)$ and $D=(\left(N-\frac{1}{2}\right)\tau,0)$ is shown in Fig. \ref{Fig:example222} (left).

Denote $\widetilde{p}:=\frac{\pi}{2}-\alpha-\varepsilon$ with a sufficiently small $\varepsilon>0$. After simple computation,
we can see that the equation of this biggest ellipse is
\begin{small}\begin{align*}
\frac{4(x-\frac{N\tau}{2})^2}{\left(\sqrt{\frac{\tau^2}{4}+\widetilde{p}^2}+\sqrt{(N-\frac{1}{2})^2\tau^2+\widetilde{p}^2}\right)^2}
+\frac{4y^2}{\left(\sqrt{\frac{\tau^2}{4}+\widetilde{p}^2}+\sqrt{(N-\frac{1}{2})^2\tau^2+\widetilde{p}^2}\right)^2-(N-1)^2\tau^2}=1.
\end{align*}\end{small}
By rescaling, we can get the axis lengths of the correspondingly biggest ellipse on and inside which $\widetilde{u}$ is analytic:
\begin{small}
$L=\sqrt{(1+\frac{1}{2(N-1)})^2+\frac{\widetilde{p}^2}{(N-1)^2}\tau^2}+\sqrt{\frac{1}{4(N-1)^2}+\frac{\widetilde{p}^2}{(N-1)^2\tau^2}}>1$,
\end{small}
and
$l=\sqrt{L^2-1}$. Thus,
\begin{equation}\label{KKI}
\tilde{K}=L+l=L+\sqrt{L^2-1}.
\end{equation}
Since $\tilde{K}>1$, we can see from (\ref{Eq:K}) that the rate of the convergence of the interpolants to $\widehat{u}(z)$ as
$ n\rightarrow \infty$ is remarkably fast.

Given $t_0>0$ and $\Lambda>1$. Now we are ready to demonstrate the estimation that introduced by the interpolation. In fact, for $0< t\leq (\Lambda-1) t_0$ and $0<\varrho<1$, together with (\ref{eq:computeu}), (\ref{eq:fullcomputeu12}), (\ref{Eq:K}) and Corollary \ref{co:delta-solution}, it is not difficult to get
\begin{align*}
\left\|u_h^N(t)-u_{I,h}^{N,n}(t)\right\|_{L^{2}(\Omega)}
\leq& C_1 L\left(\mu t_0\sin(\alpha-\tilde{d})\right)\varphi(\alpha,\tilde{d})\cdot\left(\epsilon\cdot(\epsilon_N(\varrho))^{\varrho-1}+\frac{(\epsilon_N(\varrho))^{\varrho}}{1-\epsilon_N(\varrho)}\right)\cdot \tilde{K}^{-n},
\end{align*}
where the parameters are given as mentioned in above sections.

Based on the above discussions, we can show the error estimation for the fully discrete scheme after interpolation without proof.
\begin{theorem}\label{thm:error-interpolation}
Let $u(t)$ be the solution to Problem (\ref{eq:problem}) and $u^{N,n}_{I,h}(t)$ be given by Eq. (\ref{eq:fullcomputeu12}), respectively. Given $t_0>0$ and $\Lambda>1$. For $t_0\leq t\leq (\Lambda-1) t_0$, and $0<\varrho<1$, there uniformly holds
\begin{align*}
&\left\|u(t)-u^{N,n}_{I,h}(t)\right\|_{L^{2}(\Omega)}\\
\leq & \|e_h^{N}\|_{L^2(\Omega)}+C_1\left(L(\mu t_0\sin(\alpha-\tilde{d}))\varphi(\alpha,\tilde{d})\left(\epsilon(\epsilon_N(\varrho))^{\varrho-1}
     +\frac{(\epsilon_N(\varrho))^{\varrho}}{1-\epsilon_N(\varrho)}\right)\right)\cdot \tilde{K}^{-n},
\end{align*}
where $e_h^{N}(t)$ is given in Theorem \ref{thm:homog} or Theorem \ref{Tm:fulldiscrete2}, and the parameters are given as mentioned in above sections.
\end{theorem}

Here we note that Theorem \ref{thm:error-interpolation} can tell us how many interpolation point $n$ should one choose  theoretically on a target accuracy. In the next section, we can also observe that in application $n\ll N$ is enough.

In addition, we will use the following interpolation approximation rate (IAR):
\begin{equation}\label{IAR}
IAR=\left\|u(t)-u_{I,h}^{N,n}(t)\right\|_{L^{2}(\Omega)}/\left\|u(t)\right\|_{L^{2}(\Omega)},
\end{equation}
to verify the effectiveness of our accelerated numerical scheme in the experiments section, where $u(t)$ is the analytic solution.

\begin{algorithm}[th]
    \renewcommand{\algorithmicrequire}{\textbf{Input:}}
	\renewcommand{\algorithmicensure}{\textbf{Output:}}
    \caption{CIM-Int-FEM}
		\label{CIM-FEM}
		\begin{algorithmic}[1]
           \Require $\beta$; $t_0$;
              $t_1=\Lambda t_0$; $t$; $\Lambda>1$; $\alpha$; $\delta'$; $\epsilon$; $N$; $M$; $D$.
            \Ensure the numerical result $u^{N}_{I,h}(t)$.\\
            \textbf{Subroutine 1}: Compute the optimal parameters $\varrho^*$, $\tau^*$ and $\mu^*$.
            \State \textbf{for} $j=0$; $j<D$; $j++$ $\textbf{do}$
            \State \ \ \ \ \ \ $\varrho_j$ $\leftarrow$ $1/D\cdot j;$
            \State \ \ \ \ \ \ $a(\varrho_j)$ $\leftarrow$ $\cosh^{-1}\left(\frac{\Lambda}{(1-\varrho_j)\sin(\alpha-\widetilde{d})}\right)$;
            \State \ \ \ \ \ \ $\epsilon_N(\varrho_j)$ $\leftarrow$ $\exp(\frac{2\pi \widetilde{d} N}{a(\varrho_j)})$;
            \State \ \ \ \ \ \ $\varrho^*$ $\leftarrow$
            $\min\limits_j\left(\epsilon(\epsilon_N(\varrho_j))^{\varrho_j-1}+
            \frac{(\epsilon_N(\varrho_j))^{\varrho_j}}{1-\epsilon_N(\varrho_j)}\right)$;
            \State \textbf{end do}
            \State \ \ \ \ \ \ $\tau^*$, $\mu^*$ $\leftarrow$ Substitute $\varrho^*$ into Eq. (\ref{para:Fernandz}) to  get the optimal parameters $\tau^*$ and $\mu^*$;\\
            \textbf{Subroutine 2}: Compute  the interpolation functions $\{\widehat{u}_h({z}_{x_j})\}_{j=0}^{n-1}$ on the Chebyshev points $\{z(x_j)\}_{j=0}^{n-1}$.
			\State \textbf{for} $j=0$; $j<n$; $j++$ $\textbf{do}$
            \State \ \ \ \ \ \ $x_j$ $\leftarrow$ $\frac{a+b}{2}+ \frac{b-a}{2}\cos(\frac{j\pi}{n})$ with $a=\frac{\tau*}{2}$ and $b=(N-\frac{1}{2})\tau*$;
            \State \ \ \ \ \ \ ${z}_{x_j}$ $\leftarrow$ $\mu^*(1+\sin(i x_j-\alpha))$,~ where $i^2=-1$;
            \State \ \ \ \ \ \ $\widehat{u}_h({z}_{x_j})$ $\leftarrow$ Solve the linear system: \\ \quad\quad\quad\quad\quad $\left(\left(K{z}_{x_j}+{z}^{\beta}_{x_j}\right)\widehat{u}_h({z}_{x_j}),v_h\right)+\left(\nabla \widehat{u}_{h}({z}_{x_j}), \nabla v_{h}\right)=\left(\left(K+{z}^{\beta-1}_{x_j}\right)u_{0}+\widehat{f}({z}_{x_j}),v_h\right)$;
            \State \textbf{end do}\\
            \textbf{Subroutine 3}: Compute the numerical solution $u^{N}_{I,h}(t)$.
            \State \textbf{for} $k=0$; $j<N-1$; $k++$ $\textbf{do}$
            \State \ \ \ \ \ \ $z_k$ $\leftarrow$ $\mu^*(1+\sin(i(k+1/2)\tau^*-\alpha))$,  where $i^2=-1$;
            \State \ \ \ \ \ \ $z'_k$ $\leftarrow$ $i\mu^*\cos(i(k+1/2)\tau^*-{\alpha}))$,  where $i^2=-1$;
            \State \ \ \ \ \ \ $\widehat{u}_{I,h}(z_k)$ $\leftarrow$ $\frac{\sum_{j=0}^{n}\frac{\omega_j}{z_k-{z}_{x_j}}\widehat{u}_h({z}_{x_j})}
{\sum_{j=0}^{n}\frac{\omega_j}{z_k-{z}_{x_j}}}$;
            \State \textbf{end do}
            \State \ \ \ \ \ \ $u^{N}_{I,h}(t)$ $\leftarrow$ $\frac{\tau^*}{\pi}\mathrm{Im}\left\{\sum_{k =
                            0}^{N-1}e^{z_{k}t}\widehat{u}_{I,h}(z_k)z'_{k}\right\}$.
		\end{algorithmic}
\end{algorithm}

\section{Numerical experiments}
\label{sec:Numerical}
In this section, numerical performances of our developed numerical scheme are given. During the experiments, we always
take $\Omega=(0,1)^{d}$ with $d=1,2$, $T=1$, $K=1$, $\alpha=0.6767$, $\delta'=0.1023$, and $\epsilon=2.22\times10^{-16}$.  For the spatial finite element discretization, the degree of polynomial $\Pi_{r-1}$ is taken as $r=2$.

If the exact solution $u(t)$ is known,
we define $Error_\tau(N)$ to measure the temporal errors if $h$ is given, and $Error_h(t)$ to measure the spatial errors at time $t$ if $N$ is fixed:
\begin{displaymath}
  \mathrm{Error}_{\tau}(N):=\max\limits_{t_{0}\leq t\leq \Lambda t_{0}}\left\|u(t)-u_{h}^{N}(t)\right\|_{L^{2}(\Omega)},~ ~ \mathrm{Error}_{h}(t):=\left \|u(t)-u_{h}^{N}(t)\right\|_{L^{2}(\Omega)}.
\end{displaymath}

While if the exact solution is unknown, we measure the temporal and spatial errors by using the same notations as above but defined as follows
\begin{displaymath}
  \mathrm{Error}_{\tau}(N):=\max\limits_{t_{0}\leq t\leq \Lambda t_{0}}\left\|u_{h}^{200}(t)-u_{h}^{N}(t)\right\|_{L^{2}(\Omega)},
~\mathrm{Error}_{h}(t)=\left \|u_{h}^{N}(t)-u_{h/2}^{N}(t)\right\|_{L^{2}(\Omega)}.
\end{displaymath}
Furthermore, we also use
$\frac{\ln\left(\mathrm{Error}_{h}/\mathrm{Error}_{h/2}\right)}{\ln(2)}$
to measure the convergence order in space direction.

All of the numerical experiments are implemented by using Matlab 2018a on a PC with Intel(R) Core(TM) i7-7700 CPU @3.60GHz and 4.00GB RAM.

\subsection{Example 1 (A scalar problem)}
\label{example1}
Here we aim to verify the spectral accuracy and high efficiency of CIM with the optimal parameters $\mu$, $\alpha$ and the step-spacing $\tau$ been provided in Eq. (\ref{para:Fernandz}). Firstly, we consider the scalar problem as follows:
\begin{equation}\label{Ex:example1}
K\frac{\partial u(t)}{\partial t}+\frac{1}{\Gamma(1-\beta)}\frac{\partial }{\partial t}\int_{0}^{t}(t-\tau)^{-\beta}(u(\tau)-u(0))d\tau+au(t)=f(t)~ ~ {\rm with}~ t>0,~ ~ u(0)=u_{0},
\end{equation}
and
\begin{equation*}
f(t)=1+\frac{3\sqrt{\pi}}{2}K+\frac{3\sqrt{\pi}}{2\Gamma(2-\beta)}t^{1-\beta}+\frac{3\sqrt{\pi}}{2}t,
\end{equation*}
$a=1$, $u_0=1$. The exact solution is $u(t)=1+\frac{3\sqrt{\pi}}{2}t$.

\begin{figure}[ht]
\setlength{\abovecaptionskip}{0.1cm}  
\setlength{\belowcaptionskip}{-0.2cm} 
\centering
\begin{minipage}[c]{0.39\textwidth}
 \centering
 \centerline{\includegraphics[height=3.5cm,width=4.6cm]{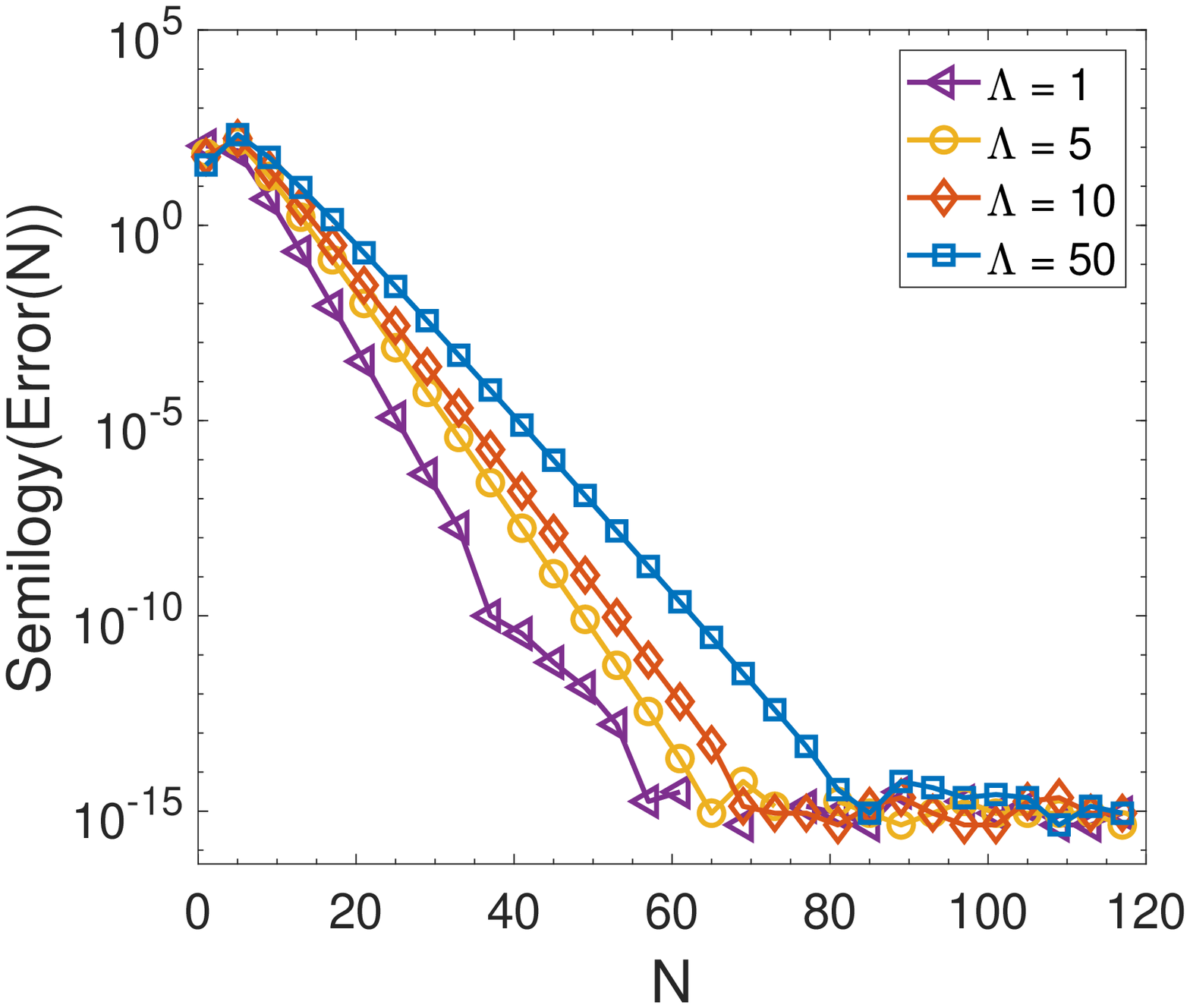}}
\end{minipage}
\begin{minipage}[c]{0.39\textwidth}
 \centering
 \centerline{\includegraphics[height=3.5cm,width=4.6cm]{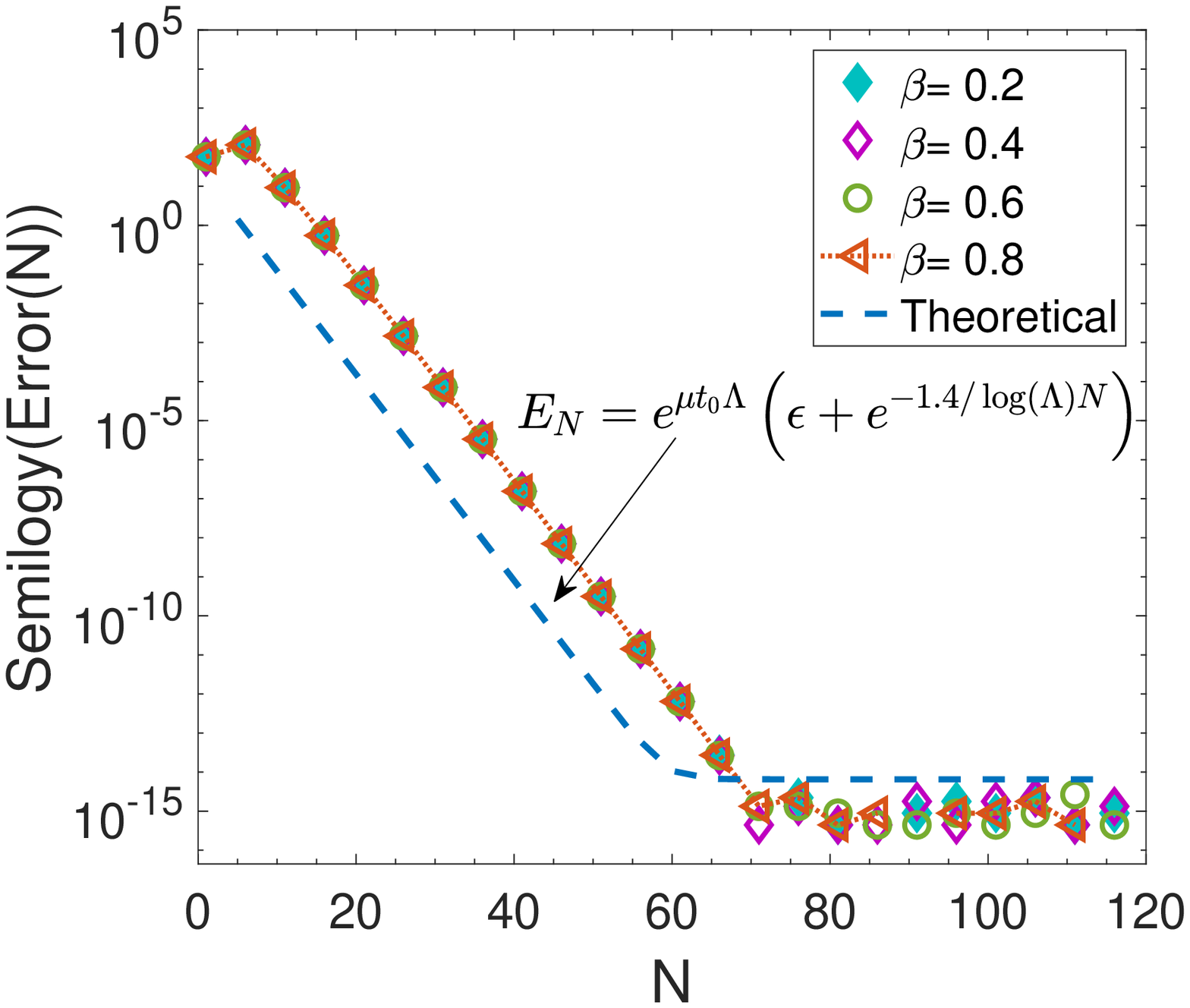}}
\end{minipage}
\caption{The numerical performance of CIM  for Problem (\ref{Ex:example1}) at time $t=0.6$.  \emph{Left: The absolute numerical errors for different $\Lambda$ with $\beta=0.5$. Right: The absolute numerical errors for different $\beta$ with $\Lambda=10$.}}
\label{Fig:example1}
\end{figure}

Define the absolute error function w.r.t. $N$ as
\begin{displaymath}
  \mathrm{Error}(N):=\max\limits_{t_{0}\leq t\leq \Lambda t_{0}}\left|u(t)-u^{N}(t)\right|.
\end{displaymath}
The absolute errors of CIM for Problem (\ref{Ex:example1}) with different parameters $\Lambda$ and fractional orders $\beta$ are shown in Fig. \ref{Fig:example1}.

One can see from Fig. \ref{Fig:example1} that the numerical results are consistent with the theoretical analysis in Theorem \ref{thm:fernandz} and Theorem \ref{EQ:convergence}. That is, CIM has spectral accuracy, and the convergence order is $\mathcal{O}(\epsilon+e^{-cN})$ with $c=\mathcal{O}(1/\log(\Lambda))$. What is more, CIM is also unconditionally stable w.r.t. the fractional order $\beta$.

\subsection{Example 2 (A problem with vanishing initial data)}
\label{sec:example2}
\begin{figure}[ht]
\centering
\begin{minipage}[c]{0.36\textwidth}
 \centering
 \centerline{\includegraphics[height=2.8cm,width=4.2cm]{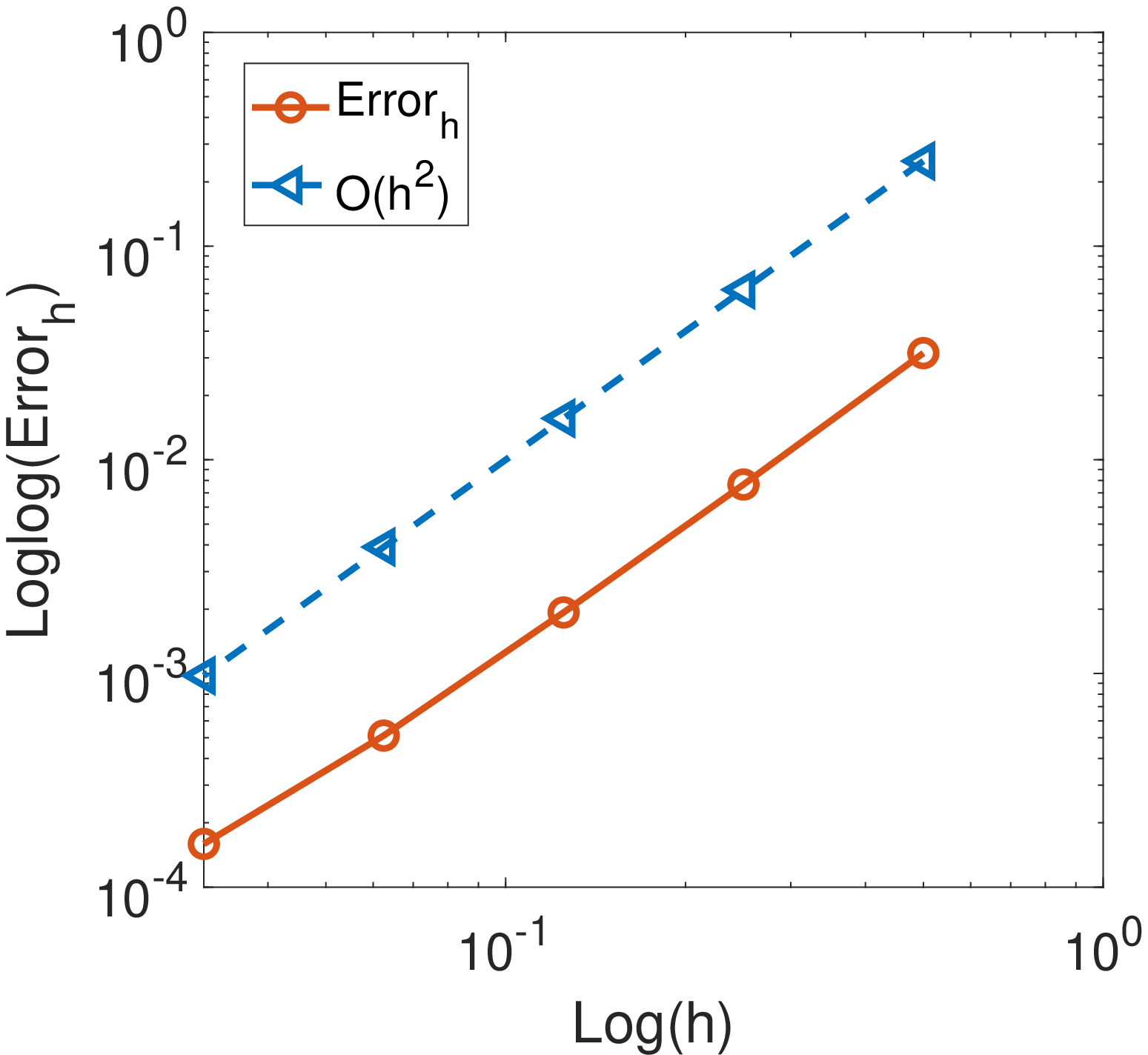}}
\end{minipage}
\begin{minipage}[c]{0.36\textwidth}
 \centering
 \centerline{\includegraphics[height=2.8cm,width=4.2cm]{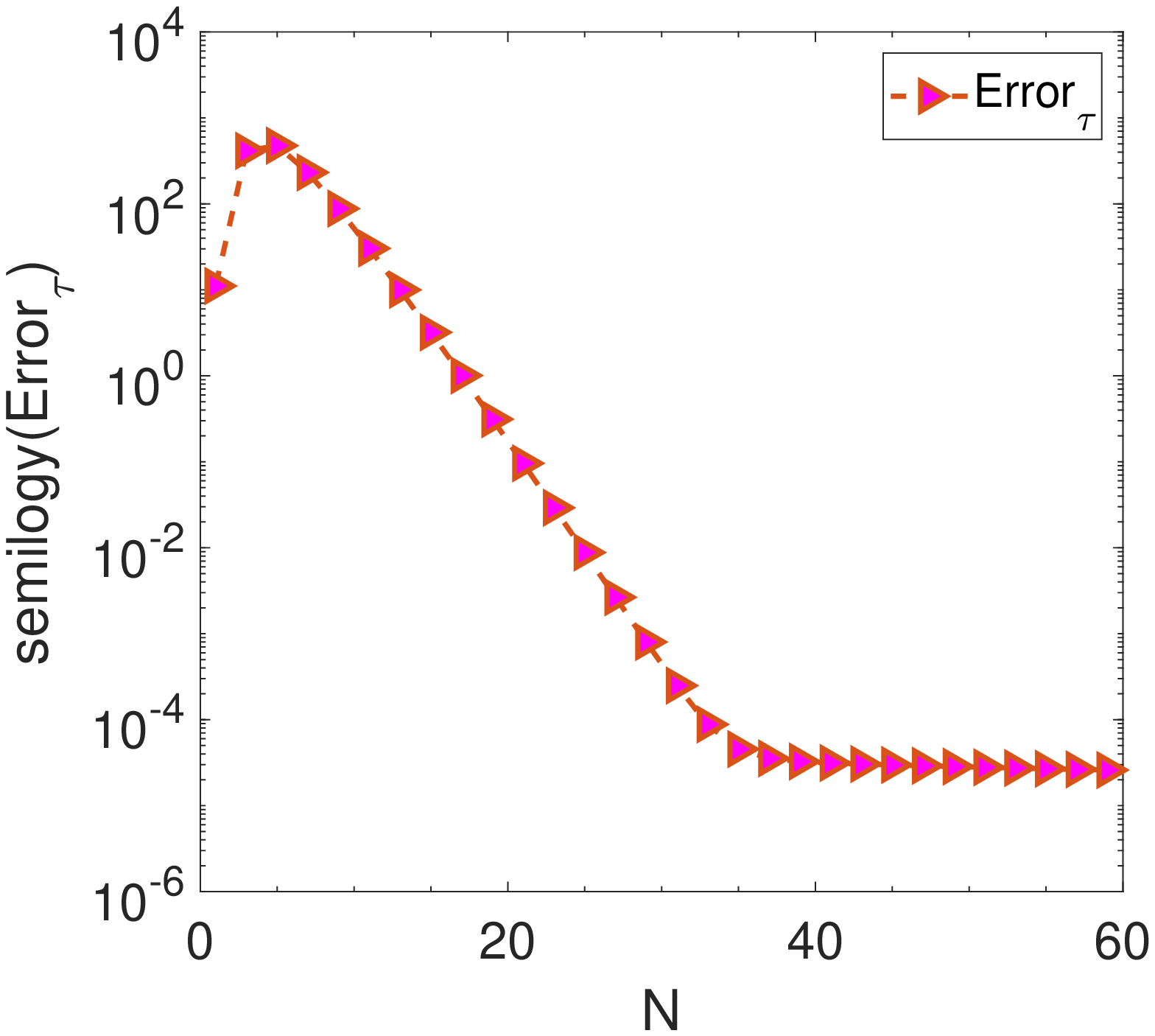}}
\end{minipage}
\begin{minipage}[c]{0.36\textwidth}
 \centering
 \centerline{\includegraphics[height=2.8cm,width=4.2cm]{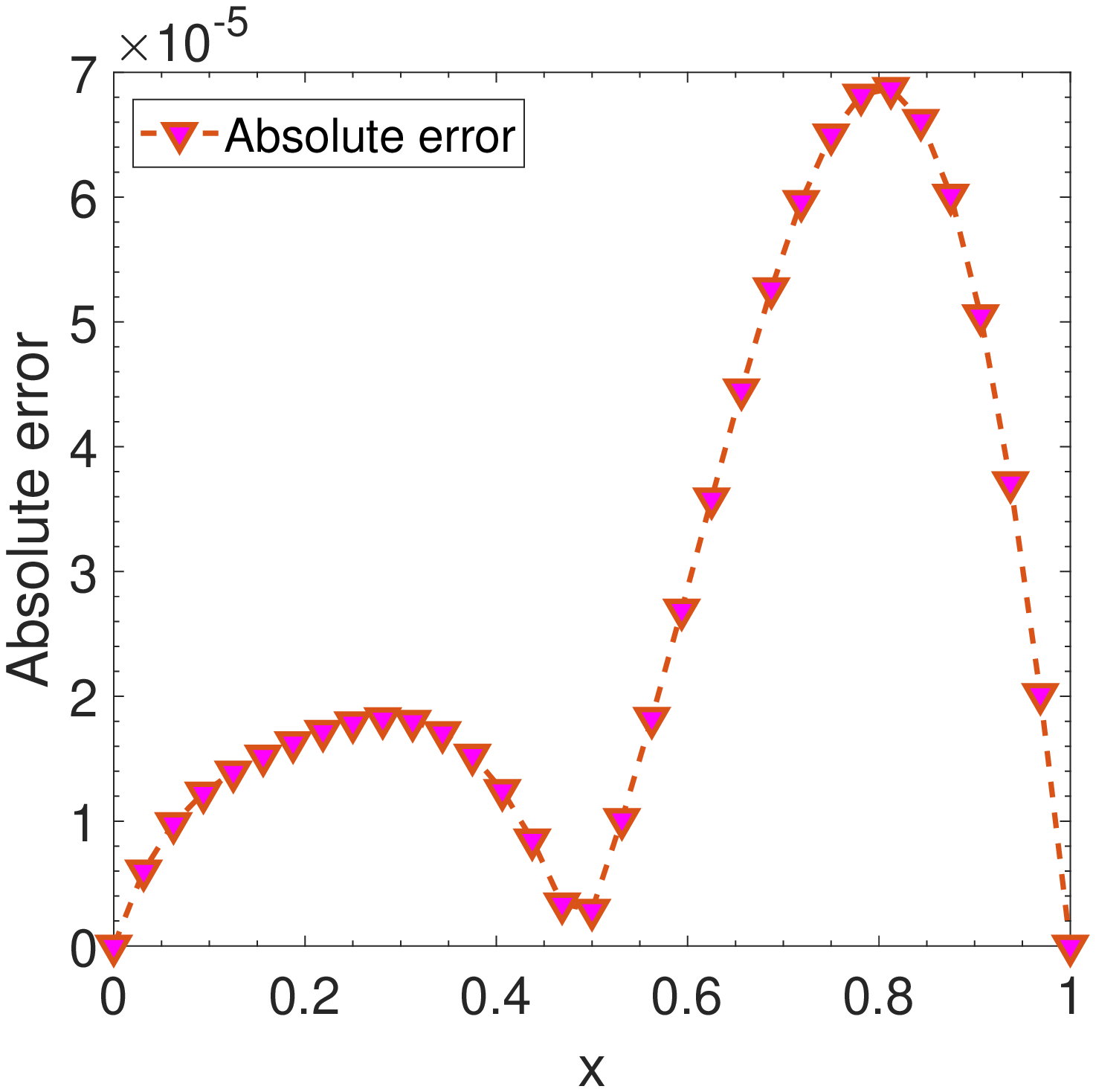}}
\end{minipage}
\begin{minipage}[c]{0.36\textwidth}
 \centering
 \centerline{\includegraphics[height=2.8cm,width=4.6cm]{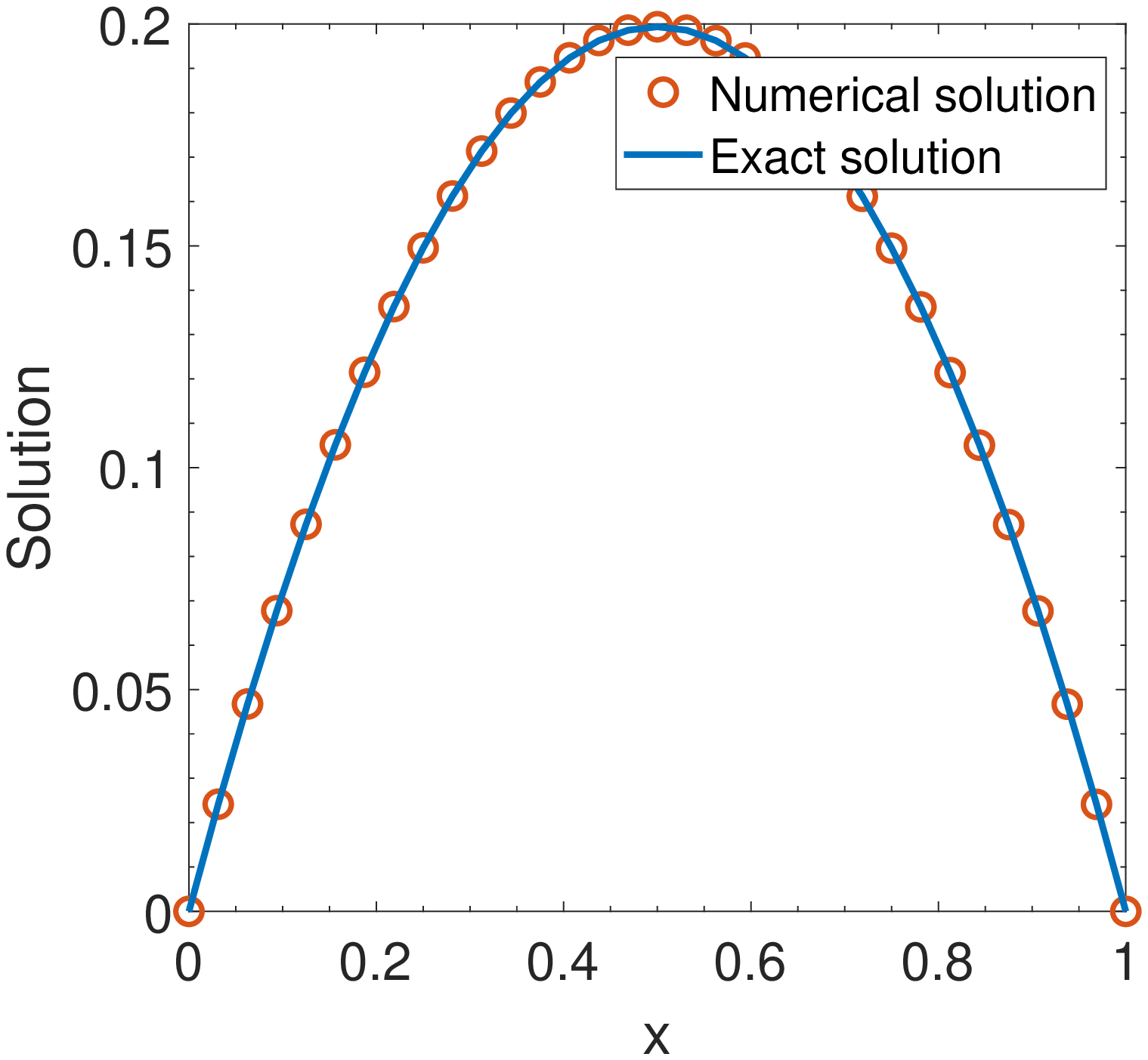}}
\end{minipage}
\caption{The numerical performance of CIM-FEM for Example \ref{sec:example2} with vanishing initial value with $\beta=0.75$ at time $t=0.86$.  \emph{(upper-left) The spatial error $\mathrm{Error}_{h}$ with $N=60$ and $h=2^{-J}$,  $J=1,2,3,4,5$. $(upper-right)$ The temporal error $\mathrm{Error}_{\tau}$ for different $N$ with $h=2^{-5}$. $(bottom-left)$ The absolute error (i.e. $|u(x,t)-u_h^{N}|$) with $h=2^{-5}$ and $N=60$. $(bottom-right)$ The fitting effect of the numerical solution with $h=2^{-5}$ and $N=60$.}}
\label{Fig:example2}
\end{figure}

To show the effectiveness of our algorithm, we consider a special case that the initial value of the problem (\ref{eq:problem}) is vanishing, and the exact solution is
\begin{displaymath}
 u(x,t)=t^{3/2}x(1-x),
\end{displaymath}
the source term is
\begin{displaymath}
f(x,t)=\frac{3}{2}Kt^{1/2}x(1-x)+\frac{3\sqrt{\pi}(5/2-\beta)}{4\Gamma(7/2-\beta)}t^{3/2-\beta}x(1-x)+2t^{3/2}.
\end{displaymath}
In this case, by choosing $\beta=0.75$ and $\Lambda=10$, the numerical performances of CIM-FEM are visually demonstrated in Fig. \ref{Fig:example2}.

We can see from Fig. \ref{Fig:example2} that convergence rate of CIM-FEM is $O(h^2)$ in space and spectral accuracy in time. In other words, the numerical results are consistent with the theoretical analysis given in Theorem \ref{thm:fernandz}, Theorem \ref{th:spaceerrorinhomo} and Theorem \ref{Tm:fulldiscrete2}.

\subsection{Example 3 (Two problems in one space dimension)}
\label{example3}
Here, we consider the 1-D homogeneous case of Problem (\ref{eq:problem}) with non-smooth as well as smooth initial data, i.e., we choose the initial data as
\begin{flushleft}
\quad\quad $(7.3.1)\quad u_0=\pi^3\chi_{(0,3/4]}(x)\in L^{2}(\Omega)$,\\
\quad\quad $(7.3.2)\quad u_0=x\cdot\chi_{[0,3/4]}(x)-x\cdot\chi_{(3/4,1]}(x)\in H^{1}_0(\Omega)$,\\
\quad\quad $(7.3.3)\quad u_0=\pi^3x(1-x)\in H^{1}_{0}(\Omega)\cap H^{2}(\Omega)$,
\end{flushleft}
respectively. To illustrate the accuracy of CIM, we take the numerical solution with $N=200$ and $h=1/2^7$ as the reference solution. The numerical results are showing in Tables \ref{Tab:01}-\ref{Tab:06}.

\begin{table*}[t!]
\setlength{\abovecaptionskip}{0.0cm}  
\setlength{\belowcaptionskip}{-0.2cm} 
\scriptsize
\caption{The temporal errors $Error_\tau$ for the homogeneous case of 1-D problem based on the non-smooth initial data (7.3.1) with $h=1/2^7$ and $t=0.8$.}
\label{Tab:01}
\begin{center}
\begin{tabular}{c|c|c|c|c|c}
\hline
\hline
\diagbox{$\beta$}{$Error_{\tau}$}{$N$} & $20$     & $40$      &$60$      &$80$       & $100$ \\
\hline
$\beta=0.25$         &2.2886E-02   & 9.8500E-05   & 9.5236E-09  & 3.8939E-14   & 3.4148E-14 \\
\hline
$\beta=0.50$         &1.9995E-02   & 9.0973E-05   & 8.8491E-09  & 8.4586E-14   & 5.4752E-14 \\
\hline
$\beta=0.75$         &1.5650E-02   & 9.0822E-05   & 7.0644E-09  & 6.7551E-14   & 6.9440E-14 \\
\hline
\hline
\end{tabular}
\end{center}
\end{table*}

\begin{table*}[t!]
\setlength{\abovecaptionskip}{0cm}  
\setlength{\belowcaptionskip}{-0.2cm} 
\scriptsize
\caption{The spatial errors $Error_h$ for the homogeneous case of 1-D problem based on the non-smooth initial data (7.3.1) with $N=100$ and $t=0.8$.}
\label{Tab:02}
\begin{center}
\begin{tabular}{c|c|c|c|c|c|c}
\hline
\hline
\multirow{2}*{$1/M$}
&\multicolumn{2}{c|}{$\beta=0.25$}&\multicolumn{2}{c|}{$\beta=0.50$} &\multicolumn{2}{c}{$\beta=0.75$}\\\cline{2-7}

              &$Error_h$     &$Order$        &$Error_h$     &$Order$     &$Error_h$    &$Order$  \\
\hline
$1/2^{5}$     & 6.3098E-04   &$--$          &1.2092E-03    &$--$       & 1.9769E-03  &$--$    \\
\hline
$1/2^{6}$     & 1.5809E-04   & 1.9969       &3.0281E-04    & 1.9976    & 4.9514E-04  &1.9974  \\
\hline
$1/2^{7}$     & 3.9543E-05   & 1.9992       &7.5735E-05    & 1.9994    & 1.2384E-04  &1.9993  \\
\hline
$1/2^{8}$     & 9.8872E-06   & 1.9998       &1.8936E-05    & 1.9998    & 3.0964E-05  &1.9998  \\
\hline
\hline
\end{tabular}
\end{center}
\end{table*}

\begin{table*}[t!]
\setlength{\abovecaptionskip}{0cm}  
\setlength{\belowcaptionskip}{-0.2cm} 
\scriptsize
\caption{The temporal errors $Error_\tau$ for the homogeneous case of 1-D problem based on the smooth initial data (7.3.2) with $h=1/2^7$ and $t=0.8$.}
\label{Tab:03}
\begin{center}
\begin{tabular}{c|c|c|c|c|c}
\hline
\hline
\diagbox{$\beta$}{$Error_{\tau}$}{$N$} & $20$     & $40$      &$60$      &$80$       & $100$ \\
\hline
$\beta=0.25$      & 1.2295E-03  & 1.4837E-03   & 1.3276E-03  & 1.0082E-03   & 8.4551E-04 \\
\hline
$\beta=0.50$      & 1.2889E-03  & 1.5165E-03   & 1.3568E-03  & 1.0299E-03   & 8.6348E-04 \\
\hline
$\beta=0.75$      & 1.3364E-03  & 1.5360E-03   & 1.3745E-03  & 1.0437E-03   & 8.7518E-04 \\
\hline
\hline
\end{tabular}
\end{center}
\end{table*}

\begin{table*}[t!]
\setlength{\abovecaptionskip}{0cm}  
\setlength{\belowcaptionskip}{-0.2cm} 
\scriptsize
\caption{The spatial errors $Error_h$ for the homogeneous case of 1-D problem based on the non-smooth initial data (7.3.1) with $N=100$ and $t=0.8$.}
\label{Tab:04}
\begin{center}
\begin{tabular}{c|c|c|c|c|c|c}
\hline
\hline
\multirow{2}*{$1/M$}
&\multicolumn{2}{c|}{$\beta=0.25$}&\multicolumn{2}{c|}{$\beta=0.50$} &\multicolumn{2}{c}{$\beta=0.75$}\\\cline{2-7}

              &$Error_h$     &$Order$        &$Error_h$     &$Order$     &$Error_h$    &$Order$  \\
\hline
$1/2^{5}$     & 6.1204E-06   &$--$          &1.1704E-05    &$--$       & 1.9077E-05   &$--$    \\
\hline
$1/2^{6}$     & 1.5355E-06   & 1.9949       &2.9347E-06    & 1.9956    & 4.7833E-06   & 1.9958  \\
\hline
$1/2^{7}$     & 3.8421E-07   & 1.9987       &7.3424E-07    & 1.9989    & 1.1967E-06   & 1.9990  \\
\hline
$1/2^{8}$     & 9.6073E-08   & 1.9997       &1.8359E-07    & 1.9997    & 2.9923E-07   & 1.9997  \\
\hline
\hline
\end{tabular}
\end{center}
\end{table*}

\begin{table*}[t!]
\setlength{\abovecaptionskip}{0cm}  
\setlength{\belowcaptionskip}{-0.2cm} 
\scriptsize
\caption{The temporal errors $Error_\tau$ for the homogeneous case of 1-D problem based on the smooth initial data (7.3.3) with $h=1/2^7$ and $t=0.8$.}
\label{Tab:05}
\begin{center}
\begin{tabular}{c|c|c|c|c|c}
\hline
\hline
\diagbox{$\beta$}{$Error_{\tau}$}{$N$} & $20$     & $40$      &$60$      &$80$       & $100$ \\
\hline
$\beta=0.25$      & 5.3750E-03  & 2.3301E-05   & 2.1672E-09  & 6.1354E-14   & 9.3454E-15 \\
\hline
$\beta=0.50$      & 4.6838E-03  & 2.1532E-05   & 2.0118E-09  & 2.4407E-14   & 1.2485E-14 \\
\hline
$\beta=0.75$      & 3.6395E-03  & 2.1530E-05   & 1.5907E-09  & 1.4396E-14   & 1.5089E-14 \\
\hline
\hline
\end{tabular}
\end{center}
\end{table*}

\begin{table*}[t!]
\setlength{\abovecaptionskip}{0cm}  
\setlength{\belowcaptionskip}{-0.2cm} 
\scriptsize
\caption{The spatial errors $Error_h$ for the homogeneous case of 1-D problem based on the smooth initial data (7.3.3) with $N=100$ and $t=0.8$.}
\label{Tab:06}
\begin{center}
\begin{tabular}{c|c|c|c|c|c|c}
\hline
\hline
\multirow{2}*{$1/M$}
&\multicolumn{2}{c|}{$\beta=0.25$}&\multicolumn{2}{c|}{$\beta=0.50$} &\multicolumn{2}{c}{$\beta=0.75$}\\\cline{2-7}

              &$Error_h$     &$Order$        &$Error_h$     &$Order$     &$Error_h$    &$Order$  \\
\hline
$1/2^{5}$     & 1.4943E-04   &$--$          &2.8646E-04    &$--$       & 4.6867E-04   &$--$    \\
\hline
$1/2^{6}$     & 3.7484E-05   & 1.9951       &7.1818E-05    & 1.9959    & 1.1749E-04   & 1.9960  \\
\hline
$1/2^{7}$     & 9.3789E-06   & 1.9988       &1.7967E-05    & 1.9990    & 2.9393E-05   & 1.9990  \\
\hline
$1/2^{8}$     & 2.3452E-06   & 1.9997       &4.4926E-06    & 1.9997    & 7.3495E-06   & 1.9998  \\
\hline
\hline
\end{tabular}
\end{center}
\end{table*}

More specifically, Tables \ref{Tab:01}, \ref{Tab:03} and  \ref{Tab:05} demonstrate the temporal errors $Error_\tau$ with non-smooth, intermediate smooth and smooth initial data separately; Tables \ref{Tab:02}, \ref{Tab:04} and Table \ref{Tab:06} list the corresponding spatial errors $Error_h$. These numerical results agree well with the theoretical analyses in Theorem \ref{thm:fernandz} with $f=0$, Theorem \ref{th:spaceerrorhomo1}, and Theorem \ref{thm:homog}.

\subsection{Example 4 (Three problems in two space dimension)}
\label{example4}
In this example, we consider the 2-D homogeneous cases of Problem (\ref{eq:problem}) with non-smooth as well as smooth initial values, i.e., we choose the initial data as
\begin{flushleft}
\quad\quad $(7.4.1)\quad u_0=\pi\chi_{(0,3/4]\times(0,1)}(x,y)\in L^{2}(\Omega)$,\\
\quad\quad $(7.4.2)\quad u_0=4\pi^2xy(1-x)(1-y)\in H^{1}_{0}(\Omega)\cap H^{2}(\Omega)$,
\end{flushleft}
respectively. In addition, we also consider the inhomogeneous problem with vanishing initial data, that is
\begin{flushleft}
\quad\quad $(7.4.3)\quad f(x,y,t)=3\pi^5 e^{\frac{3}{2}t}\sin(x)(1-x)^{2}y(y-1)$, and $u_0=0$.
\end{flushleft}
The numerical results are presented in Tables \ref{Tab:07}-\ref{Tab:12}. Among them, Tables \ref{Tab:07}, \ref{Tab:09} and \ref{Tab:11} show the temporal numerical errors $Error_\tau$ for the three cases above; Tables \ref{Tab:08}, \ref{Tab:10} and \ref{Tab:12} illustrate the corresponding spatial errors $Error_h$. These numerical results verify Theorem \ref{thm:fernandz}, Theorem \ref{th:spaceerrorinhomo}, Theorem \ref{thm:homog} and Theorem \ref{Tm:fulldiscrete2}.
\begin{table*}[t!]
\setlength{\abovecaptionskip}{0cm}  
\setlength{\belowcaptionskip}{-0.2cm} 
\scriptsize
\caption{The temporal errors $Error_\tau$ for the homogeneous case of 2-D problem based on the non-smooth initial data (7.4.1) with $h=1/2^5$ and $t=0.6$.}
\label{Tab:07}
\begin{center}
\begin{tabular}{c|c|c|c|c|c}
\hline
\hline
\diagbox{$\beta$}{$Error_{\tau}$}{$N$} & $20$     & $40$      &$60$      &$80$       & $100$ \\
\hline
$\beta=0.25$         &2.3321E-03  & 1.4808E-03   & 3.7972E-04  & 6.0484E-05    & 8.1762E-06 \\
\hline
$\beta=0.50$         &3.2892E-03  & 1.6862E-03   & 3.0516E-04  & 5.3018E-05    & 7.8812E-06 \\
\hline
$\beta=0.75$         &4.2886E-03  & 1.6947E-03   & 2.4295E-04  & 4.1240E-05    & 6.5707E-06 \\
\hline
\hline
\end{tabular}
\end{center}
\end{table*}
\begin{table*}[t!]
\setlength{\abovecaptionskip}{0cm}  
\setlength{\belowcaptionskip}{-0.2cm} 
\scriptsize
\caption{The spatial errors $Error_h$ for the homogeneous case of 2-D problem based on the non-smooth initial data (7.4.1) with $N=100$ and $t=0.6$.}
\label{Tab:08}
\begin{center}
\begin{tabular}{c|c|c|c|c|c|c}
\hline
\hline
\multirow{2}*{1/M}
&\multicolumn{2}{c|}{$\beta=0.25$}&\multicolumn{2}{c|}{$\beta=0.50$} &\multicolumn{2}{c}{$\beta=0.75$}\\\cline{2-7}
              &$Error_h$    &$Order$   &$Error_h$     &$Order$   &$Error_h$    &$Order$\\
\hline
$1/2^{3}$     &1.6503E-03   &$--$     &1.8865E-03    &$--$     &2.0332E-03   &$--$   \\
\hline
$1/2^{4}$     &3.9126E-04   &2.0765   &4.7358E-04    &1.9941   &5.5716E-04   &1.8676 \\
\hline
$1/2^{5}$     &9.6537E-05   &2.0190   &1.1877E-04    &1.9955   &1.4253E-04   &1.9668 \\
\hline
$1/2^{6}$     &2.4044E-05   &2.0054   &2.9700E-05    &1.9996   &3.5821E-05   &1.9924 \\
\hline
\hline
\end{tabular}
\end{center}
\end{table*}

\begin{table*}[t!]
\setlength{\abovecaptionskip}{0cm}  
\setlength{\belowcaptionskip}{-0.2cm} 
\scriptsize
\caption{The temporal errors $Error_\tau$ for the homogeneous case of 2-D problem based on the smooth initial data (7.4.2) with $h=1/2^5$ and $t=0.6$.}
\label{Tab:09}
\begin{center}
\begin{tabular}{c|c|c|c|c|c}
\hline
\hline
\diagbox{$\beta$}{$Error_{\tau}$}{$N$} & $20$     & $40$      &$60$      &$80$       & $100$ \\
\hline
$\beta=0.25$         &1.1401E-03  & 6.8451E-04   & 1.6312E-04  & 2.5065E-05   & 3.4184E-06 \\
\hline
$\beta=0.50$         &1.6228E-03  & 7.6423E-04   & 1.2958E-04  & 2.2021E-05   & 3.2967E-06 \\
\hline
$\beta=0.75$         &2.1277E-03  & 7.4689E-04   & 1.0162E-04  & 1.7137E-05   & 2.7505E-06 \\
\hline
\hline
\end{tabular}
\end{center}
\end{table*}

\begin{table*}[t!]
\setlength{\abovecaptionskip}{0cm}  
\setlength{\belowcaptionskip}{-0.2cm} 
\scriptsize
\caption{The spatial errors $Error_h$ for the homogeneous case of 2-D problem based on the  smooth initial data (7.4.2) with $N=100$ and $t=0.6$.}
\label{Tab:10}
\begin{center}
\begin{tabular}{c|c|c|c|c|c|c}
\hline
\hline
\multirow{2}*{1/M}
&\multicolumn{2}{c|}{$\beta=0.25$}&\multicolumn{2}{c|}{$\beta=0.50$} &\multicolumn{2}{c}{$\beta=0.75$}\\\cline{2-7}
              &$Error_h$    &$Order$   &$Error_h$    &$Order$   &$Error_h$    &$Order$\\
\hline
$1/2^{3}$     & 7.2485E-04  &$--$     & 7.4908E-04   &$--$     & 6.8250E-04   &$--$   \\
\hline
$1/2^{4}$     & 1.5771E-04  &2.2004   & 1.5307E-04   &2.2909   & 1.2654E-04   &2.4312 \\
\hline
$1/2^{5}$     & 3.7754E-05  &2.0626   & 3.5949E-05   &2.0902   & 2.8830E-05   &2.1340 \\
\hline
$1/2^{6}$     & 9.3307E-06  &2.0166   & 8.8396E-06   &2.0239   & 7.0319E-06   &2.0356 \\
\hline
\hline
\end{tabular}
\end{center}
\end{table*}

\begin{table*}[t!]
\setlength{\abovecaptionskip}{0cm}  
\setlength{\belowcaptionskip}{-0.2cm} 
\scriptsize
\caption{The temporal errors $Error_\tau$ for the inhomogeneous case of 2-D problem based on the source term (7.4.3) with $u_0=0$, $h=1/2^5$ and $t=0.6$.}
\label{Tab:11}
\begin{center}
\begin{tabular}{c|c|c|c|c|c}
\hline
\hline
\diagbox{$\beta$}{$Error_{\tau}$}{$N$} & $20$     & $40$      &$60$      &$80$       & $100$ \\
\hline
$\beta=0.25$         &7.5385E-02  & 4.3031E-02   & 2.2798E-02  & 8.9194E-03   & 3.2602E-03 \\
\hline
$\beta=0.50$         &9.4276E-02  & 4.7014E-02   & 2.0328E-02  & 8.3595E-03   & 3.2009E-03 \\
\hline
$\beta=0.75$         &1.1184E-01  & 4.7663E-02   & 1.8007E-02  & 7.3730E-03   & 2.9225E-03 \\
\hline
\hline
\end{tabular}
\end{center}
\end{table*}

\begin{table*}[t!]
\setlength{\abovecaptionskip}{0cm}  
\setlength{\belowcaptionskip}{-0.2cm} 
\scriptsize
\caption{The spatial errors $Error_h$ for the inhomogeneous case of 2-D problem based on the source term (7.4.3) with $u_0=0$, $N=100$ and $t=0.6$.}
\label{Tab:12}
\begin{center}
\begin{tabular}{c|c|c|c|c|c|c}
\hline
\hline
\multirow{2}*{$1/M$}
&\multicolumn{2}{c|}{$\beta=0.25$}&\multicolumn{2}{c|}{$\beta=0.50$} &\multicolumn{2}{c}{$\beta=0.75$}\\\cline{2-7}
              &$Error_h$    &$Order$   &$Error_h$   &$Order$   &$Error_h$    &$Order$\\
\hline
$1/2^{3}$     & 1.0311E-02  &$--$     &1.4646E-02   &$--$     & 1.8830E-02   &$--$   \\
\hline
$1/2^{4}$     & 2.6138E-03  &1.9799   &3.7230E-03   &1.9760   & 4.7974E-03   &1.9727 \\
\hline
$1/2^{5}$     & 6.5389E-04  &1.9990   &9.3200E-04   &1.9980   & 1.2016E-03   &1.9973 \\
\hline
$1/2^{6}$     & 1.6346E-04  &2.0001   &2.3303E-04   &1.9998   & 3.0048E-04   &1.9996 \\
\hline
\hline
\end{tabular}
\end{center}
\end{table*}

\subsection{Example 5 (Numerical performance of the acceleration algorithm)}
This experiment targets on verifying the efficient performance of the acceleration algorithm developed in Section \ref{sec:acceler}. For simplicity, ``{CIM-FEM}'' indicates that the designed numerical scheme is implemented directly without interpolation; ``{CIM-Int-FEM}'' represents the algorithm equipped with the barycentric interpolation (see Section \ref{sec:acceler}).

\begin{figure}[ht]
\setlength{\abovecaptionskip}{0cm}  
\setlength{\belowcaptionskip}{-0.2cm} 
  \centering
  \includegraphics[height=7.5cm,width=10cm]{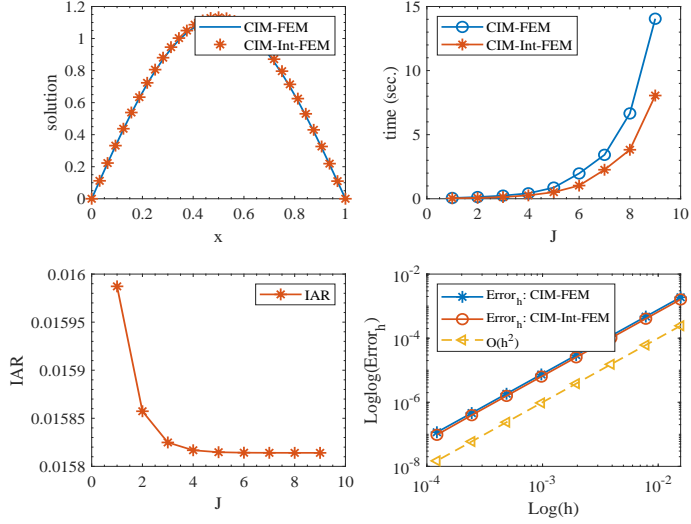}\\
\caption{Numerical comparisons of the CIM with and without interpolation in 1-D situation by solving Problem (\ref{eq:problem}) with initial data (7.3.1) when taking $N=100$. (\emph{upper-left}) describes the fitting effect of the two numerical solutions when $h=2^{-5}$; (\emph{upper-right}) depicts the CPU time costs of CIM-Int-FEM and CIM-FEM for different $h$; (\emph{bottom-left}) shows the trend of IAR by using CIM-Int-FEM as $h$ decreases; (\emph{bottom-right}) represents the spatial convergence order of the CIM with and without interpolation.}
  \label{Fig:performance}
\end{figure}

\begin{figure}[ht]
\setlength{\abovecaptionskip}{0cm}  
\setlength{\belowcaptionskip}{-0.2cm} 
  \centering
  \includegraphics[height=7.5cm,width=10cm]{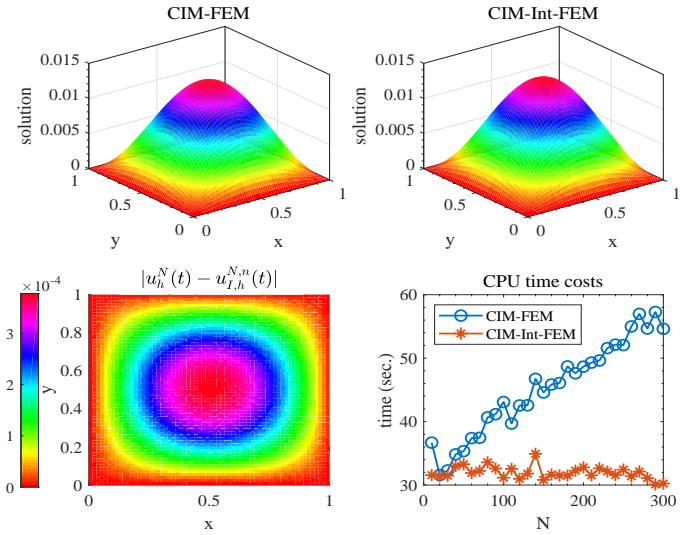}\\
\caption{Numerical performance of the CIM with and without interpolation in 2-D situation by solving Problem (\ref{eq:problem}) with initial data (7.4.2) when taking $h=2^{-6}$. (\emph{upper-left}) and (\emph{upper-right}) are the numerical solutions $u_h^{N}(t)$ computed by CIM-FEM and $u_{I,h}^{N,n}(t)$ by CIM-Int-FEM respectively with $N=100$, ; (\emph{bottom-left}) shows the error $|u_h^{N}(t)-u_{I,h}^{N,n}(t)|$; (\emph{bottom-right}) represents the CPU time costs of CIM-Int-FEM and CIM-FEM for different $N$.}
\label{Fig:performance2D}
\end{figure}

We compare these two algorithms by solving Problem (\ref{eq:problem}) with both of the non-smooth initial value (7.3.1) for 1-D case and the smooth initial value (7.4.2) for 2-D case, respectively. The number of the interpolation nodes is chosen as $n=10$, and $\beta=0.5$, $t=0.4$. In the 1-D case, the spatial step-spacing is $h=2^{-(J+4)}$, $J\in\{1,2,...,9\}$.
Detailed comparisons are shown in Fig. \ref{Fig:performance} and Fig. \ref{Fig:performance2D}.

Fig. \ref{Fig:performance} and Fig. \ref{Fig:performance2D} indicate the effectiveness of interpolation in computational time. What is more, we can see from the bottom-right subplot of Fig. \ref{Fig:performance} that the CIM after acceleration does not effect the spatial accuracy, and from the bottom-right one of Fig. \ref{Fig:performance2D} that the rate of convergence in $N$ is not deteriorated by interpolation.

In short, all of the above examples verify that both of the schemes CIM-Int-FEM and CIM-FEM we proposed are effective. The spatial finite element discretization can reach to the optimal second-order convergence and the time CIM discretization has spectral accuracy. The improved acceleration interpolation process is significant.


\section{Conclusions}
\label{sec:conclusion}
In this paper, the time fractional normal-subdiffusion transport model, which depicts a crossover from
normal diffusion (as $t\rightarrow 0$) to sub-diffusion (as $t\rightarrow \infty$), is analyzed and numerically solved.
First, based on the analytic properties of the bivariate Mittag-Leffler function, we rigorously discuss the regularity, existence and uniqueness  of the model solution of the homogeneous as well as inhomogeneous cases. In particular, new regularity results show that, with time developing, the time regularity of the solution to TFDEs is dominated by the highest-order operator.
Then, a numerical scheme with high-precision and low regularity requirements is developed. More specifically, for time discretization, by Laplace transform, a parallel CIM scheme is designed, in which the hyperbolic integral contour is selected; for space discretization, the standard Galerkin finite element method is employed. Error estimates show that CIM-FEM has spectral accuracy in time and can reach to optimal second-order convergence with smooth as well as non-smooth initial data in space. Furthermore, the barycentric interpolation algorithm is proposed to accelerate the algorithm, and the error estimate shows that the convergence rate of this interpolation is remarkably fast. Several numerical experiments in 1-D and 2-D for both homogeneous and inhomogeneous cases verify all of the theoretical results. High numerical performances of these numerical examples also illustrate the high efficiency and robustness of our numerical schemes.
\begin{flushleft}
\textbf{Data Availability} Data sharing not applicable to this article as no datasets were generated or analysed during the current study.
\end{flushleft}

\section*{Acknowledgements}
This research was supported by National Natural Science Foundation of China under Grant Nos. 12071195 and 12225107, and the Innovative Groups of Basic Research
in Gansu Province under Grant No. 22JR5RA391. The author Zhao is supported by
Guangdong Basic and Applied Basic Research Foundation No. 2022A1515011332 and the Fundamental Research Funds for the Central University under Grant No. D5000230096. The authors have no relevant financial or non-financial interests to disclose.

\section*{Appendix \textbf{A}: Limiting asymptotic behavior of the MSD of subordinator $\mathbf{B}[\mathbb{E}_t]$}

Let $\mathbf{B}(t)$ denote the $d$-dimensional Brownian motion and $\mathbb{E}_t$ be the inverse process or
hitting time process of the $\beta$-stable subordinator $\mathbb{S}_t= Kt +\overline{S}_t$ with drift, where $0<\beta<1$ and $K>0$ ($K=0$ implies drift-less). Then, by the celebrated L\'{e}vy-Khintchine formula \cite{Applebaum09}, we have
$\mathbb{E}[e^{i(\lambda,\mathbb{S}_t)}]=e^{-t(K\lambda+\lambda^{\beta})}$.

Denote $G$ and $H$ as the PDFs of $\mathbf{B}(t)$ and $\mathbb{E}_t$, respectively. Then, there is $u(t,x)=\int_{0}^{\infty}G(x,s)H(s,t)ds$. Since $\partial G(x,t)=\Delta G(x,t)$ and $H(u,t)=-\partial_u\int_{0}^{t}g(y,u)dy$ with $g$ being the PDF of $\mathbb{S}_t$, then by Laplace transformation, we obtain $\widehat{H}(u,z)=-\partial_u\widehat{g}(u,z)z^{-1}=\frac{Kz+z^{\beta}}{z}e^{-u(Kz+z^{\beta})}$. Based on these facts, the mean-squared displacement (MSD) of $\mathbf{B}[\mathbb{E}_t]$ is
\begin{displaymath}\begin{split}
D(t):=\left\langle(\mathbf{B}[\mathbb{E}_t])^2\right\rangle=\int_{\mathbb{R}^d}|x|^{2}u(t,x)dx
=\int_{0}^{\infty}H(u,t)\int_{\mathbb{R}^d}|x|^{2}G(x,u)dxdu=2d\langle\mathbb{E}_t\rangle,
\end{split}\end{displaymath}
where $d$ is the dimension of space. Transforming $\langle\mathbb{E}_t\rangle$ by Laplace w.r.t. $t$, we have
\begin{displaymath}
\mathcal{L}_{t\rightarrow z}\{\langle\mathbb{E}_t\rangle\}(z):=\int_{0}^{\infty}u\widehat{H}(u,z)du=\frac{1}{z(Kz+z^{\beta})}\simeq
\left\{\begin{split}
&K^{-1}z^{-2},\quad && z\rightarrow\infty,\\
&z^{-(1+\beta)},\quad &&z\rightarrow 0.
\end{split}\right.
\end{displaymath}
Further, by Tauberian's theorem (cf. \cite[Theorem 1.5.7]{Applebaum09}), the limiting asymptotic behaviors of MSD are
\begin{displaymath}
D(t)\simeq
\left\{\begin{split}
&2dK^{-1}t,\quad &&t\rightarrow 0,\\
&\frac{2d}{\Gamma(1+\beta)}t^{\beta},&&t\rightarrow\infty.
\end{split}\right.
\end{displaymath}
Now, it can be clearly seen that the model we discuss in this paper depicts an crossover from normal diffusion to sub-diffusion. That is, when $t$ is small enough, it portrays normal diffusion; as $t$ grows large enough, it reflects sub-diffusion.

\section*{Appendix \textbf{B}}

\subsection*{{\rm\textbf{B.1 Definition of the multivariate Mittag-Leffler function}}}
\label{sebsec:111}

Let $\alpha,\beta,\gamma \in\mathbb{C}$, and ${\rm Re}(\alpha)>0,{\rm Re}(\beta)>0$. The bivariate Mittag-Leffler function is defined as
\begin{equation}\label{def:BIML}
E_{(\alpha,\beta),\gamma}^{\delta}(z_1, z_2):=\sum\limits_{k=0}^{\infty}\sum\limits_{l=0}^{\infty}\frac{(\delta)_{k+l}}{\Gamma(\alpha k+\beta l+\gamma)}\cdot\frac{z_1^k z_2^l}{k!l!},
\end{equation}
where the numerator $(\delta)_{k+l}$ is the Pochhammer symbol, i.e.,
\begin{displaymath}
(a)_n=\frac{\Gamma(a+n)}{\Gamma(a)}=a(a+1)(a+2)...(a+n-1).
\end{displaymath}
As mentioned in \cite{Bazhlekova21}, the multiple power series (\ref{def:BIML}), which converges absolutely and locally uniformly, defines an entire function in $z_1$ and $z_2$.

\subsection*{\textbf{B.2 Proof of Lemma \ref{Lem:MandL}}}
\begin{proof}
Let $\alpha$, $\beta$, $\gamma$, $\omega_1$, $\omega_2\in \mathbb{R}$ with $0<\alpha<\beta\leq1$, $\gamma\geq0$ and $\omega_1$, $\omega_2<0$. Firstly, for $t=0$, according to the series representation of the bivariate Mittag-Leffler function in Eq. (\ref{def:BIML}), there exist a positive constant $C$ such that
\begin{equation}\label{est:M_LL}
\left|E^1_{(\alpha,~\beta),~\gamma}\left(\omega_1 t^{\alpha}, \omega_2 t^{\beta}\right)\right|=\frac{1}{\Gamma(\gamma)}\leq C.
\end{equation}

For $t>0$, by Corollary 2 in \cite{Fernandez20}, the contour integral representation of the bivariate Mittag-Leffler function $E^1_{(\alpha,~\beta),~\gamma}(\omega_1 t^{\alpha}, \omega_2 t^{\beta})$ is
\begin{equation}\begin{split}\label{est:M_LLL}
E^1_{(\alpha,~\beta),~\gamma}\left(\omega_1 t^{\alpha}, \omega_2 t^{\beta}\right)
=\frac{t^{1-\gamma}}{2\pi i}\int_{\Gamma_{\theta,\delta}}\frac{e^{zt}z^{-\gamma}}{1+|\omega_1|z^{-\alpha}+|\omega_2|z^{-\beta}}dz.
\end{split}\end{equation}
\begin{description}
  \item[Case I] for given $\beta\in(0,1)$ and any $\alpha\in(0,1)$ with $0<\alpha<\beta<1$. When $z\in\Gamma_{\theta,\delta}$, we choose $\theta\in(\frac{\pi}{2},\pi)$, which only depends on $\beta$ and closes to $\frac{\pi}{2}$ enough such that $\arg(|\omega_2|z^{-\beta})>-\frac{\pi}{2}$. Then the angle between $|\omega_2|z^{-\beta}$ and $1+|\omega_1|z^{-\alpha}$ is less then $\pi/2$ and the denominator in (\ref{est:M_LLL}) satisfies
  \begin{equation}\label{est:alphaa}
      \left|1+|\omega_1|z^{-\alpha}+|\omega_2|z^{-\beta}\right|\geq|\omega_2||z|^{-\beta}\cos\beta\theta.
  \end{equation}
 Based on these, choosing $\delta=1/t>0$ large enough, then we get
\begin{displaymath}\begin{split}
 \left|E^1_{(\alpha,~\beta),~\gamma}\left(\omega_1 t^{\alpha}, \omega_2 t^{\beta}\right)\right|
&\leq \frac{t^{1-\gamma}}{ 2\pi}\int_{\Gamma_{\theta,\delta}}\frac{e^{|z|t}|z|^{-\gamma}}{\left|1+|\omega_1|z^{-\alpha}+|\omega_2|z^{-\beta}\right|}|dz|\\
&\leq \frac{t^{1-\gamma}}{2\pi\left|\omega_2\right|\cos(\beta\theta)}
\int_{\Gamma_{\theta,\delta}}e^{|z|t}|z|^{\beta-\gamma}|dz| \\
&\leq\frac{t^{1-\gamma}}{2\pi\left|\omega_2\right|\cos(\beta\theta)}Ct^{\gamma-1-\beta}
\leq\frac{C}{|\omega_1t^{\beta}|}.
\end{split}\end{displaymath}
  \item[Case II] for given $\alpha\in(0,1)$, $\beta=1$. Similar to (i), when $z\in\Gamma_{\theta,\delta}$, by choosing $\theta\in(\frac{\pi}{2},\pi)$, which depends on $\alpha$ and closes to $\frac{\pi}{2}$ enough such that $\alpha\theta+\arg(1+|\omega_2|z^{-1})\geq -\frac{\pi}{2}$.
  In this case, the denominator in (\ref{est:M_LLL}) satisfies
  \begin{equation}\label{est:alphab}
      \left|1+|\omega_2|z^{-1}+|\omega_1|z^{-\alpha}\right|>\left|1+|\omega_2|z^{-1}\right|
      \geq|\omega_2||z|^{-1}\cos\left(\theta-\frac{\pi}{2}\right).
  \end{equation}
  Similarly, there holds
\begin{displaymath}\begin{split}
 \left|E^1_{(\alpha,~\beta),~\gamma}\left(\omega_1 t^{\alpha}, \omega_2 t^{\beta}\right)\right|
&\leq \frac{t^{1-\gamma}}{ 2\pi}\int_{\Gamma_{\theta,\delta}}\frac{e^{|z|t}|z|^{-\gamma}}{\left|1+|\omega_1|z^{-\alpha}+|\omega_2|z^{-\beta}\right|}|dz|\\
&\leq \frac{t^{1-\gamma}}{2\pi\left|\omega_2\right|\sin(\theta)}
\int_{\Gamma_{\theta,\delta}}e^{|z|t}|z|^{1-\gamma}|dz|\\
&\leq\frac{t^{1-\gamma}}{2\pi\left|\omega_2\right|\sin(\theta)}Ct^{\gamma-1-1}
\leq\frac{C}{|\omega_2t|}.
\end{split}\end{displaymath}
\end{description}
Thus the proof of this lemma is completed. \qed
\end{proof}

\section*{Appendix \textbf{C}: Proof of Corollary \ref{co:delta-solution}}
\begin{proof}
Firstly, we let $t_0<t\leq \Lambda t_0$ and take the solution to Problem (\ref{eq:problem}) as $u(x,t)=\delta(t-t_0)\cdot\mathbf{1}_{(\Omega)}$, where $\delta(t)$ is the Dirac delta function. Then $u_0\equiv 0$, $\widehat{u}(x,z)=e^{-zt_0}\cdot\mathbf{1}_{(\Omega)}$, $\widehat{f}(x,z)=(Kz+z^{\beta})e^{-zt_0}\cdot\mathbf{1}_{(\Omega)}$, and $$u^{N}(t)=\frac{\tau}{2\pi i}\sum_{k=-(N-1)}^{N-1}e^{z(\phi_k)t}e^{-z(\phi_k)t_0}z'(\phi_k).$$

From (\ref{norm_f})) and (\ref{ieq:inequality1}), there is a constant $C_2$, which depends on $\Omega,K,\beta,\alpha,\mu$, such that
$$\left\|\widehat{f}(z)\right\|_{Nr}
\leq\sup\limits_{z\in D\subset\Sigma_{\theta}}|\Omega|\left|Kz+z^{\beta}\right|\left|e^{-zt_0}\right|
\leq C|z|\cdot\left|e^{-zt_0}\right|\leq C_2.$$
Thus, for $t_0< t\leq\Lambda t_0$, by Theorem \ref{thm:fernandz}, we can obtain
\begin{align*}
    &\left\|\frac{\tau}{2\pi i}\sum_{k=-(N-1)}^{N-1}e^{z(\phi_k)(t-t_0)}z'(\phi_k)\right\|_{L^2(\Omega)}\\
\leq& C_2 L\left(\mu t_0\sin(\alpha-\tilde{d})\right)\varphi(\alpha,\tilde{d})\cdot\left(\epsilon\cdot(\epsilon_N(\varrho))^{\varrho-1}+\frac{(\epsilon_N(\varrho))^{\varrho}}{1-\epsilon_N(\varrho)}\right),
\end{align*}
i.e.,
\begin{align*}
    &\left\|\sum_{k=-(N-1)}^{N-1}e^{z(\phi_k)(t-t_0)}z'(\phi_k)\right\|_{L^2(\Omega)}\\
\leq & C_1 L\left(\mu t_0\sin(\alpha-\tilde{d})\right)\varphi(\alpha,\tilde{d})\cdot\left(\epsilon\cdot(\epsilon_N(\varrho))^{\varrho-1}+\frac{(\epsilon_N(\varrho))^{\varrho}}{1-\epsilon_N(\varrho)}\right).
\end{align*}
Actually, the above estimation is obtained by using the fact that $\int_{0}^{\infty}e^{-\mu t\sin(\alpha-\tilde{d})\cosh(x)}\leq  \int_{0}^{\infty}e^{-\mu t_0\sin(\alpha-\tilde{d})\cosh(x)}\leq L(\mu t_0\sin(\alpha-\tilde{d}))$ (cf. Lemma \ref{et:lemma}). In other words,
\begin{align*}
    &\left\|\sum_{k=-(N-1)}^{N-1}e^{z(\phi_k)(t-t_0)}z'(\phi_k)\right\|_{L^2(\Omega)}\\
\leq & C_1\int_{0}^{\infty}e^{-\mu t\sin(\alpha-\tilde{d})\cosh(x)}\varphi(\alpha,\tilde{d})\cdot\left(\epsilon\cdot(\epsilon_N(\varrho))^{\varrho-1}+\frac{(\epsilon_N(\varrho))^{\varrho}}{1-\epsilon_N(\varrho)}\right),
\end{align*}
still holds.

So, for $0< t\leq (\Lambda-1) t_0$, there is
\begin{align*}
    &\left\|\sum_{k=-(N-1)}^{N-1}e^{z(\phi_k)t}z'(\phi_k)\right\|_{L^2(\Omega)}\\
\leq & C_1\int_{0}^{\infty}e^{-\mu (t+t_0)\sin(\alpha-\tilde{d})\cosh(x)}\varphi(\alpha,\tilde{d})\cdot\left(\epsilon\cdot(\epsilon_N(\varrho))^{\varrho-1}+\frac{(\epsilon_N(\varrho))^{\varrho}}{1-\epsilon_N(\varrho)}\right)\\
\leq & C_1L\left(\mu t_0\sin(\alpha-\tilde{d})\right)\varphi(\alpha,\tilde{d})\cdot\left(\epsilon\cdot(\epsilon_N(\varrho))^{\varrho-1}+\frac{(\epsilon_N(\varrho))^{\varrho}}{1-\epsilon_N(\varrho)}\right).
\end{align*}
\end{proof}

%
%


\end{document}